\RequirePackage{etoolbox}
\csdef{input@path}{{style/}{graphics/}}
\documentclass[aap,MSNbibl,seceqn,dvips]{arximspdf}
\makeatletter
   \@ifpackageloaded{graphicx}{}{\usepackage{graphicx}}
\makeatother

\doi{10.1214/14-AAP1029} %kopijuoti is PTS
\volume{25}
\issue{3}
\pubyear{2015}
\firstpage{1513}
\lastpage{1539}

\makeatletter
\newcommand{\eqref}[1]{(\ref{#1})}
\newtheorem{theorem}{Theorem}[section]
\newtheorem{prop}[theorem]{Proposition}
\newtheorem{cor}[theorem]{Corollary}
\newtheorem{lemma}[theorem]{Lemma}
\newproclaim{assumption}{Assumption}
\newproclaim{remark}{Remark}
\newproclaim{example}{Example}

\def\mean{\mathbb{E}}
\def\Reals{\mathbb{R}}
\def\Ints{\mathbb{Z}}
\def\Nats{\mathbb{N}}
\def\Coms{\mathbb{C}}

\def\tailr{r} % Tail ratio
\def\ldh{H} % First ladder height
\def\iof{K} % Integrate overflow
\def\renf{U} % Renewal measure and function
\def\cdf{F} % Cumulative distribution function
\def\rvf{A} % Regularly varying function
\def\svf{\beta} % Slowly varying function
\def\CR{\mathcal{R}}
\def\iff{\Longleftrightarrow}

\def\rx{\varepsilon}

\def\dpois{\mathrm{Poisson}}
\def\levy{\mbox{L\'evy}}

\def\iunit{\mathrm{i}}

\makeatother

\begin{document}
\begin{frontmatter}

\title{Strong renewal theorems with infinite mean beyond~local~large~deviations}
\runtitle{SRT beyond LLD}

\begin{aug}
\author[A]{\fnms{Zhiyi}~\snm{Chi}\corref{}\ead[label=e1]{zhiyi.chi@uconn.edu}}
\runauthor{Z. Chi}
\affiliation{University of Connecticut}
\address[A]{Department of Statistics\\
University of Connecticut\\
215 Glenbrook Road, U-4120\\
Storrs, Connecticut 06269\\
USA\\
\printead{e1}}
\end{aug}

% HISTORY:
\received{\smonth{6} \syear{2013}}
\revised{\smonth{3} \syear{2014}}

% ABSTRACT
%
\begin{abstract}
Let $\cdf$ be a distribution function on the line in the domain of
attraction of a stable law with exponent $\alpha\in(0,1/2]$. We
establish the strong renewal theorem for a random walk $S_1, S_2,
\ldots $
with step distribution $\cdf$, by extending the large deviations
approach in Doney [\textit{Probab. Theory Related Fileds} \textbf{107} (1997)
451--465]. This is done by introducing
conditions on $\cdf$ that in general rule out local large
deviations bounds of the type $\mathbb{P}\{S_n \in(x,x+h]\} = O(n)
\overline{\cdf}(x)/x$, hence are significantly weaker than the
boundedness condition in Doney (1997). We also give
applications of the results on ladder height processes and
infinitely divisible distributions.
\end{abstract}
%

% KEYWORDS

\end{frontmatter}

%s1 #&#
\section{Introduction} \label{sintro}
Let $X,X_1, X_2, \ldots $ be i.i.d. real-valued random variables with
distribution function $\cdf$. Denote $S_n=\sum_{i=1}^n X_i$. This
article concerns the asymptotic of
%
%e1.1 #&#
\begin{equation}
\label{erenew-measure} \renf(x+I):=\sum_{n=1}^\infty
\mathbb{P}\{ S_n \in x+I\} \qquad\mbox{as } x \to\infty,
\end{equation}
under certain tail conditions on $\cdf$, where $I = (0,h]$ with $h\in
(0,\infty)$. Specifically, denoting by $\CR_\alpha$ the class of
functions that are regularly varying at $\infty$ with exponent
$\alpha$ and $\overline{\cdf}(x) = 1-\cdf(x) = \mathbb{P}\{X>x\}$,
the first
condition is
%
%e1.2 #&#
\begin{equation}
\label{etail-svf}\overline{\cdf}(x) \sim1/\rvf(x)\qquad
 \mbox{as }x\to\infty  \mbox{ with }
\rvf \in\CR_\alpha, \alpha\in(0,1).
\end{equation}
By \eqref{etail-svf}, $p^+:= \mathbb{P}\{X>0\}>0$. The second
condition is the tail ratio condition
%
%e1.3 #&#
\begin{equation}
\label{etail-ratio} \tailr_F:=\lim_{x\to\infty } \bigl\{
\cdf(-x)/\overline{\cdf}(x)\bigr\}\qquad  \mbox{exists and is finite}.
\end{equation}
Actually, we often only need the following weaker tail ratio condition:
%
%e1.4 #&#
\begin{equation}
\label{etail-ratio-2} \limsup_{x\to\infty } \bigl\{\cdf(-x)/\overline {\cdf}
(x)\bigr\} < \tailr< \infty.
\end{equation}

There are several well-known works on the strong renewal theorem
(SRT) for~$S_n$, that is, the nontrivial limit of $x \overline{\cdf}(x)
\renf(x+I)$ as $x\to\infty $ with $0<h<\infty$; see
\cite{garsialamperti63,williamson68,doney97} for the arithmetic
case and \cite{erickson70tams} for the nonarithmetic case. The
definitions of (non)arithmetic distributions and the related
(non)lattice distributions are given in Section~\ref{sresults}. While the
SRT always holds for $\alpha\in(1/2,1]$ in the arithmetic case as
well as in the nonarithmetic case with the extra condition $\mathbb
{P}\{X\ge
0\} = 1$, there are examples where it fails to hold for $\alpha\in
(0,1/2]$; see \cite{williamson68}, and also
\cite{garsialamperti63} for explanations. For the arithmetic case,
a well-known condition that leads to the SRT for all $\alpha\in
(0,1/2]$ is
\[
\sup_{n\ge0} \bigl\{n\mathbb{P}\{X=n\}/\overline{\cdf}(n)\bigr
\}<\infty,
\]
provided $X$ is integer-valued. Under the condition, the SRT was
established for $1/4<\alpha\le1/2$ in \cite{williamson68}.
The general arithmetic case remained open until
\cite{doney97}, which took a different approach from previous
efforts. The core of the argument in \cite{doney97} is an estimate
of local large deviations (LLD) for the events $\{S_n\in x+I\}$ as
$x\to\infty $. Once these estimates are established, the rest of the proof
is basically an application of the local limit theorems (LLTs)
(\cite{bingham89}, Theorem~8.4.1--2) \cite
{garsialamperti63,williamson68,erickson70tams}. Recently, it was shown
\cite{vatutin13tpa} that, for the nonlattice case, if
%
%e1.5 #&#
\begin{equation}
\label{edoney} \sup_{x\ge0} \omega_I(x) < \infty,
\end{equation}
where for $I\subset\Reals$,
%
%e1.6 #&#
\begin{equation}
\label{eprob-ratio} \omega_I(x) = x \mathbb{P}\{X\in x+I\}/\overline {
\cdf}(x),
\end{equation}
then a much simpler argument than the one in \cite{doney97} can be
used to get the same type of LLD bound, which then leads to the SRT.

However, condition \eqref{edoney} can be restrictive. As an
example, let $\cdf$ be supported on $[1, \infty)$ and have piecewise
constant density $f(x) \propto h(x)$, such that
\[
h(x) = \cases{ n^{-\alpha-1}, &\quad $n\le x < n+1, n\in\Nats\setminus \bigl
\{2^k, k\ge1 \bigr\},$ \vspace*{2pt}
\cr
k n^{-\alpha-1}, & \quad $n =
2^k \le x < n+1,  k\in\Nats$.}
\]
Then $\overline{\cdf}(x) \sim C/x^{-\alpha}$ as $x\to\infty $,
where $C, C',
\ldots$ denote constants. Set $I=(0,1]$. Since $\omega_I(x)\sim
C'\ln x$ for $x\in[2^k, 2^k+1)$, \eqref{edoney} does not
hold. On the other hand, the set of $x$ with large $\omega_I(x)$ has
low density of order $O(\ln x/x)$, while the large values of
$\omega_I(x)$ increases slowly at order $O(\ln x)$ as $x\to\infty $. Thus
it is reasonable to wonder if the SRT should still hold.

To handle similar situations as the example, one way is to control the
aggregate effect of large values of $\omega_I(x)$. We therefore
define the function
%
%e1.7 #&#
\begin{equation}
\label{eoverflow} \iof(x,T) = \iof(x, T; I)= \int_0^x
\bigl[\omega_I(y) - T \bigr]^+\,\mathrm{d}y,
\end{equation}
where $c^\pm= \max(\pm c,0)$ for $c\in\Reals$ and $T>0$ is a
parameter. We will show that, for example, if $X>0$ and $\alpha\in
(0,1/2)$, and if for some $T>0$,
%
%e1.8 #&#
\begin{equation}
\label{eiof} \iof(x,T) = o \bigl(\rvf(x)^2 \bigr),
\end{equation}
then the SRT holds for $S_n$. Since \eqref{edoney} implies
$\iof(x,T)\equiv0$ if $T>0$ is large enough, it is a special case of
\eqref{eiof}. In the above example, since $\iof(x,T) = O((\ln x)^2)$
for large $T>0$, the SRT holds as well. Notice that if \eqref{eiof}
holds for one $h\in(0,\infty)$, it holds for all $h\in(0,\infty)$.
As will be seen, the condition can be further relaxed.

There has been constant interest in the large deviations of sums of
random variables with regularly varying distributions and infinite
mean ($0<\alpha<1$), most notably in the ``big-jump'' domain; see
\cite{hult05aap,denisov08ap,doney97,doney12ptrf,tang06ejp} and
references therein. The theme of this line of
research is to identify the domain of large $x$, such that the event
$S_n\in x+I$ with $a_n\ll x$ and $0<h\le\infty$ is mainly due to a
single large value among $X_1, \ldots, X_{n}$. Here $a_n$ are
constants such
that $S_n/a_n$ is tight; see the definition of $a_n$ in
Section~\ref{sresults}. The local version of this type of large
deviations, with $h<\infty$ as opposed to $h=\infty$ in the global
version, requires more elaborate conditions on $\mathbb{P}\{X\in x+I\}
$, and
it seems that none of the conditions in the current literature allows
occasional large values of $\omega_{I}(x)$. As shown in \cite
{doney97,vatutin13tpa}, to establish the SRT, the precise LLD
$\mathbb{P}\{S_n\in x+I\} \sim n \mathbb{P}\{X\in x + I\}$ is
unnecessary, and instead
$\mathbb{P}\{S_n\in x+I\} = O(n) \mathbb{P}\{X\in x + I\}$ or even
$\mathbb{P}\{S_n\in x+I\} =
O(n) \overline{\cdf}(x)/x$ can be the starting point. Our results
implies the latter is not necessary either; in the above example, for
each fixed large $n$, $\limsup_{x\to\infty } \mathbb{P}\{S_n\in
x+I\}/[n
\overline{\cdf}(x)/x] = \infty$, because as $k\to\infty $,
\begin{eqnarray*}
&& \mathbb{P}\bigl\{S_n\in\bigl[2^k+a_n,
2^k+a_n+1\bigr)\bigr\}
\\
&&\qquad\ge n \mathbb{P}\bigl\{X\in\bigl[2^k, 2^k+1/2\bigr)\bigr\}
\mathbb{P}\bigl\{S_{n-1} \in[a_n, a_n+1/2)\bigr
\}
\\
& &\qquad\sim c n a_n^{-1} k/2^{k(\alpha+1)}
\end{eqnarray*}
and $\overline{\cdf}(2^k+a_n)/(2^k + a_n) \sim c'/2^{k(\alpha+1)}$, where
$c, c'>0$ are constants. On the other hand, as in \cite
{doney97}, we still need certain estimates of the $\levy$
concentration function of $S_n$ \cite{petrov95}. These are
systematically furnished by the analysis on small-step sequences in
\cite{denisov08ap}.

As an application, we will consider the ladder height process of
$S_n$. Because $\cdf$ is the basic information, it is of interest to
find conditions on $\iof$ that yield the SRT for the ladder height
process. It is known that under certain conditions, the step
distribution of the ladder height process is in the domain of
attraction of stable law with exponent $\alpha\varrho$, where
$0\le\varrho\le1$ is the positivity parameter of the limiting stable
distribution associated with $S_n$ \cite{rogozin71}.
We will show that, if $\iof(x,T) = O(\rvf(x)^{2c})$ for some $c\in
[0,\varrho)$, then the SRT holds for the ladder process. Note that
since the ladder steps are nonnegative, due to the results in \cite
{garsialamperti63,erickson70tams}, only the case where $\alpha
\varrho\le1/2$ needs to be considered.

As another application, we will also consider the case where $X$ is
infinitely divisible. Since the $\levy$ measure $\nu$ of $X$ is
typically much easier to specify than its distribution function
$\cdf$, a natural question is whether similar conditions on $\iof(x,
T)$ can be found for $\nu$ that lead to the SRT. This question turns
out to have a positive answer. Naturally, it is more interesting
and important to study the SRT for $\levy$ processes under a similar
setting. However, this is beyond the scope of the paper.

The main results are stated in Section~\ref{sresults} and their proofs are
given in Section~\ref{sproof-renewal}.

%s2 #&#
\section{Main results} \label{sresults}
Since other than \eqref{etail-svf}, there are no constraints on
$\rvf(x)$, we shall always assume without loss of generality that
it is continuous and strictly increasing on $[0,\infty)$ with
$\rvf(0)>0$, such that for $x\gg1$, $\rvf(x) = x^\alpha
\exp \{\int_1^x \rx(v) \,\mathrm{d}v/v \} $, where $\rx
(v)$ is bounded and
continuous, and $\rx(v)\to0$ as $v\to\infty $ (\cite{korevaar04},
Theorem~IV2.2). Then
%
%e2.1 #&#
\begin{equation}
\label{einverse} \rvf^{-1}(x) = x^{1/\alpha} \svf(x)
\end{equation}
is continuous and strictly increasing, with $\svf\in\CR_0$
(\cite{bingham89}, Theorem~1.5.12). Furthermore, by $\rvf'(x) \sim
\alpha\rvf(x)/x$, $(\rvf^{-1})'(x) \sim\rvf^{-1}(x)/(\alpha x)$.
Denote $a_n = \rvf^{-1}(n)$. Then, as $n\to\infty $, $a_n\to\infty
$ and
$n(1-a_n/a_{n+1}) = 1/\alpha+ o(1)$.

Under conditions \eqref{etail-svf} and \eqref{etail-ratio}, $S_n /
a_n \stackrel{D}{\rightarrow} \zeta$, where $\zeta$ is a stable
random variable such
that (\cite{breiman92}, pages~207--213)
\begin{eqnarray}
\mean \bigl[e^{\iunit\theta\zeta} \bigr] &=& \exp \biggl\{\int
\bigl(e^{\iunit\theta x}-1 \bigr)\lambda(x) \,\mathrm{d}x \biggr\}\nonumber
\\
\eqntext{\mbox{with } \lambda(x) = \mathbf{1}\{x>0\} x^{-\alpha-1} +
\tailr_F \mathbf{1}\{x<0\} |x|^{-\alpha-1}.} %
\end{eqnarray}
Letting $\gamma=
(1-\tailr_F)/(1+\tailr_F)$, the positivity parameter $\varrho=
\mathbb{P}\{\zeta>0\}$ of $\zeta$ is equal to $1/2 + (\pi\alpha)^{-1}
\tan^{-1}(\gamma\tan(\pi\alpha/2))$ (\cite{bingham89}, page~380).
Let $g$ denote the density of~$\zeta$.

Henceforth, $\cdf$ is said to be arithmetic (resp., lattice), if there
is $d>0$, such that its support is contained in $d\Ints$ (resp., $a +
d\Ints$ for some $0\le a<d$). In either case, the span of $\cdf$ is
the largest such $d$. $\cdf$ is said to be nonarithmetic (resp.,
nonlattice) if it is not arithmetic (resp., not lattice). A lattice
distribution can be nonarithmetic. Indeed, provided the span of the
distribution is $d$, the distribution is nonarithmetic $\iff$ its
support is contained in $a+d\Ints$ for some $a>0$ with $a/d$ being an
irrational number.

We shall always assume $h>0$ is fixed. Since it is well-known that
for $\alpha\in(1/2,1)$, the SRT holds if (1) $\cdf$ is nonarithmetic
and concentrated in $[0,\infty)$ \cite{erickson70tams}, or (2) $\cdf$
is arithmetic \cite{williamson68,garsialamperti63}, we shall only
consider $\alpha\in(0,1/2]$.

The main results of this section are obtained under the following:

%as1 #&#
\begin{assumption} \label{abasic}
There exist a function $L$ and a constant $T_0>0$ such that,
letting $\theta= 1/\alpha-1$, the following hold:
\begin{longlist}[(a)]%[itemsep=0ex, leftmargin=0ex,
%itemindent=2.4\parindent, label={\alph*})]
%
\item[(a)] %\label{aL}
$L\in\CR_c$ for some $c\in[0,\alpha]$ and is nondecreasing. If
$p^+=1$, then $L(x)\to\infty $. If $p^+\in(0,1)$, then $L(x)/\ln
x\to\infty $. Furthermore,
%
%e2.2 #&#
\begin{equation}
\label{elow-cut} x\overline{\cdf}(x) \sum_{n\le L(x)}
\mathbb{P}\{S_n\in x+I\}\to0, \qquad x\to\infty.
\end{equation}
\item[(b)] %\label{aiof-1}
If $\alpha\in(0,1/2)$, then
%
%e2.3 #&#
\begin{equation}
\label{ediff2} \iof(x, T_0) = o \biggl( \frac{\rvf(x)^2}{u_\theta
(x)} \biggr),\qquad
\mbox{where } u_\theta(x) = \sum_{n\ge L(x)}
\frac{n^{-\theta}}{\svf(n)};
\end{equation}
\item[(c)] %\label{aiof-2}
If $\alpha=1/2$, then
%
%e2.4 #&#
\begin{eqnarray}
\label{ediff} %
\iof(x, T_0) = \cases{ \displaystyle O \biggl(
\frac{\rvf(x)^2}{\tilde u(x)} \biggr),
 &\quad $ \displaystyle\frac{\tilde u(x)}{\tilde
u(L(x))}\to1,$ \vspace*{2pt}
\cr
\displaystyle o \biggl(
\frac{\rvf(x)^2}{\tilde u(x)} \biggr), &\quad  $\mbox{else},$}
\nonumber
\\[-8pt]
\\[-8pt]
\eqntext{\mbox{where }\displaystyle  \tilde u(x) = \int_1^{\rvf(x)}
\frac{y^{-1} \,\mathrm{d}y}{\svf(y)}. }
\end{eqnarray}
\end{longlist}
\end{assumption}

%th2.1 #&#
\begin{theorem} \label{trenewal}
Let $\alpha\in(0,1/2]$ and \eqref{etail-svf}--\eqref
{etail-ratio} hold. Then Assumption~\ref{abasic} implies the SRT
%
%e2.5 #&#
\begin{equation}
\label{erenewal} %
\lim_{x\to\infty } x\overline{\cdf}(x)
\renf(x+I) = h\Lambda_\cdf\qquad
\mbox{with } \Lambda_\cdf= \alpha\int_0^\infty
x^{-\alpha} g(x) \,\mathrm{d}x,
\end{equation}
where $h>0$ is arbitrary if $\cdf$ is nonarithmetic, and is the span
of $\cdf$ otherwise.
\end{theorem}

%re1 #&#
\begin{remark*}
\begin{longlist}[(4)] %[itemsep=0ex, leftmargin=0ex,
%itemindent=2.2\parindent, label=\arabic*.]
%
\item[(1)]
If $\tailr_F=0$ in \eqref{etail-ratio}, then $\Lambda_\cdf=
\sin(\pi\alpha)/\pi$; see \cite{erickson70tams}.

\item[(2)]
Under Assumptions~\ref{abasic}(a) and (b),
$u_\theta\in\CR_{c(2-1/\alpha)}$. Since $\theta= 1/\alpha-1>1$,
if $c>0$, then clearly the order of the bound in \eqref{ediff2} is
strictly higher than $\rvf(x)^2$. If $c=0$, then by $L(x)\to\infty $,
$u_\theta(x)=o(1)$, so the bound in \eqref{ediff2} is still
strictly higher than $\rvf(x)^2$.

\item[(3)]
In Assumption~\ref{abasic}(c), $\tilde u(x)$ is
increasing in $x>0$. Also, either $\tilde u\in\CR_0$ or
$\tilde u(x)$ converges to a finite number as $x\to\infty $.

\item[(4)]
The integral conditions in \eqref{ediff2} and \eqref{ediff} also
imply some ``hard'' upper limits to $\omega_I(x)$. Indeed, since
uniformly for $t\in[0,h]$,
\begin{eqnarray*}
&&\omega_I(x-t) + \omega_I(x+h-t)
\\
&&\qquad\sim\frac{x}{\overline{\cdf}(x)} \bigl[\mathbb{P}\{X\in x-t+I\} + \mathbb {P}\{X\in x +
h-t+I\}\bigr] \ge \omega_I(x)
\end{eqnarray*}
as $x\to\infty $, if $\omega_I(x_n)\to\infty $ for a sequence
$x_n\to\infty $, then
for any $T>0$,
%
%e2.6 #&#
\begin{eqnarray}
\label{eiof-o} \iof(x_n+h, T) &\ge&\int_0^h
\bigl( \bigl[\omega_I(x_n-t) -T \bigr]^++ \bigl[
\omega_I(x_n+h-t)-T \bigr]^+ \bigr) \,\mathrm{d}t
\nonumber
\\[-8pt]
\\[-8pt]
\nonumber
&\ge& h \bigl[\omega_I(x_n)-2T \bigr]^+.
\end{eqnarray}
Therefore, the bound in \eqref{ediff2} or \eqref{ediff} applies
to $\omega_I(x)$ as well. %\qed
\end{longlist}
\end{remark*}

It can be shown that if the SRT holds, then \eqref{elow-cut} holds
for any $L(x) = o(\rvf(x))$; see the \hyperref[app]{Appendix}. The
question is,
before validating the SRT for $\cdf$, whether one can find $L$ so
that \eqref{elow-cut} holds, and if so, how fast $L$ can
grow? It is easy to see that if
%
%e2.7 #&#
\begin{equation}
\label{etail-power0} %
 x \overline{\cdf}(x)
\mathbb{P}\{X \in x + I\} \to0\quad
\mbox{or, equivalently} \quad \omega_I(x) = o \bigl(
\rvf(x)^2 \bigr),
\end{equation}
then \eqref{elow-cut} holds provided $L$ grows slowly enough, and
so Assumption~\ref{abasic}(a) is satisfied if $p^+=1$. Also
note that by \eqref{eiof-o}, $\omega_I(x) = O(\iof(x,T_0))$. Then
the next result is immediate.

%co2.2 #&#
\begin{cor} \label{crenewal}
Let $p^+=1$ and \eqref{etail-svf}--\eqref{etail-ratio} hold.
Let $\alpha\in(0,1/2)$ or $\alpha=1/2$ and $\tilde u(\infty) <
\infty$. If \eqref{etail-power0} holds and $\iof(x,T_0) =
O(\rvf(x)^2)$, in particular, if $\iof(x, T_0) = o(\rvf(x)^2)$, then
the SRT holds for $\renf$.
\end{cor}

%ex1 #&#
\begin{example*}
In \cite{williamson68}, it is shown that if $X$
only takes values in $\Nats$, such that
\[
\mathbb{P}\{X=n\} = \cases{ C n^{-3/2}\ln n, &\quad $n\ne2^k
\mbox{ for some $k\in\{0\}\cup\Nats$},$ \vspace*{2pt}
\cr
C n^{-1/2}/(\ln
n)^q, &\quad $ \mbox{otherwise, with $q=1$,}$}
\]
where $C=C(q)>0$ is a constant that may change from line to line,
then \eqref{etail-power0} does not hold, and hence the SRT fails to
hold. We show that if $q\ge2$, then the SRT holds. First, as in
\cite{williamson68}, $\overline{\cdf}(x) \sim x^{-1/2}\ln x$, where $C>0$
is a constant. Then $\alpha=1/2$, \eqref{etail-power0} holds, and we
can set $\rvf(x) = C\sqrt{x}/\ln x$. Since $X$ is aperiodic with
support $\{0\}\cup\Nats$, $h=1$. It follows that $\rvf^{-1}(x) =
x^2\svf(x) \sim C (x\ln x)^2$. By setting $T_0>0$ large enough,
$[\omega_I(x) - T_0]^+>0$ if and only if $x\in[2^k-1, 2^k)$ for
some $k\in\Nats$, and in this case, $[\omega_I(x) - T_0]^+\sim C
x/(\ln x)^{1+q}$. Then $\iof(x, T_0) \sim C x/(\ln x)^q$. Because
$\svf(x)\sim(\ln x)^2$, it is easy to check that $\tilde u(x)$
converges as $x\to\infty $. Then by Corollary~\ref{crenewal}, the SRT
holds for $\renf$. %\qed
\end{example*}

On the other hand, if $p^+\in(0,1)$, our argument for Theorem~\ref{trenewal} requires $L$ grow faster than $\ln x$. Meanwhile,
it is desirable to have faster growth of $L$ in order to get weaker
conditions on $\iof(x,T_0)$. We have the following prior lower bound
on the growth of $L$.

%
%pr2.3 #&#
\begin{prop} \label{pprior}
Let \eqref{etail-svf} and \eqref{etail-ratio-2} hold with $\alpha
\in(0,1)$. If for some $\kappa\in[0, 2\alpha)$,
%
%e2.8 #&#
\begin{equation}
\omega_I(x) = O \bigl(x^\kappa \bigr), \label{etail-power}
\end{equation}
then for any
$\rx\in(0, 2\alpha-\kappa)$, \eqref{elow-cut} holds with
$L(x) = x^{\rx/2}$.
\end{prop}

Now we consider the SRT for the ladder height processes of $S_n$, with
$S_0=0$. The (strict) ascending ladder height process $\ldh_n$ is
defined to be $S_{T_n}$, where $T_0=0$, $T_n = \min\{k\dvtx S_k>\ldh
_{n-1}\}$, $n\ge1$. The weak ascending ladder height process
is defined by replacing $>$ with $\ge$ in the definition of $T_n$, and
the descending process is defined by symmetry. Then $\ldh_n$ is a
random walk such that the steps are i.i.d. $\sim\ldh_1$. Denote
$\ldh= \ldh_1$ and $\renf^+$ the renewal measure for $\ldh_n$.
As noted in the \hyperref[sintro]{Introduction}, in the next
statement, we explicitly
require $\alpha\varrho\in(0,1/2]$.

%th2.4 #&#
\begin{theorem} \label{tladder}
Let $\alpha\in(0,1)$ and \eqref{etail-svf}--\eqref{etail-ratio}
hold. Let $\alpha\varrho\in(0,1/2]$ and
either \textup{(a)} $\varrho\in(0,1)$, or \textup{(b)} $\varrho=1$ and $S_n\to\infty
$ a.s.
If there exist $T>0$ and $c\in[0,\varrho)$, such that
%
%e2.9 #&#
\begin{equation}
\iof(x, T) = o \bigl(\rvf(x)^{2c} \bigr),
\end{equation}
then the SRT holds for $\renf^+$,
\[
\lim_{x\to\infty } x \mathbb{P}\{\ldh>x\} \renf^+(x+I) = h \sin (\pi
\alpha \varrho)/\pi,
\]
where $h>0$ is arbitrary if $\cdf$ is nonarithmetic, and is the span
of $\cdf$ otherwise.
\end{theorem}

%re2 #&#
\begin{remark*}
Under the same conditions, the SRT also holds for the weak ladder
process. %\qed
\end{remark*}

Now suppose $X$ is infinitely divisible with $\levy$ measure $\nu$,
such that
\[
\mean \bigl[e^{\iunit\theta X} \bigr] = \exp \biggl\{ \iunit\mu \theta-
\sigma^2\theta^2/2 + \int \bigl(e^{\iunit\theta u} - 1 -
\iunit\theta u\mathbf{1}\bigl\{|u|\le1\bigr\} \bigr) \nu(\mathrm{d}u) \biggr\}.
\]
For $x>0$, denote $\overline{\nu}(x) = \nu((x,\infty))$ and $\nu
(-x) =
\nu((-\infty, -x))$. Define
\[
\tilde\iof(x,T) = \int_0^x \bigl[\tilde
\omega_I(y) - T \bigr]^+\,\mathrm{d}y\qquad \mbox{where } \tilde
\omega_I(x) = x\nu(x+I)/\overline{\nu}(x).
\]

%th2.5 #&#
\begin{theorem} \label{trenewal-ind}
Let
%
%e2.10 #&#
\begin{equation}
\label{etail-ind} \overline{\nu}(x)\sim1/\rvf(x), \nu(-x) \sim \tailr
_\nu/ \rvf(x),\qquad x\to\infty,
\end{equation}
where $\rvf\in\CR_\alpha$ with $\alpha\in(0,1/2]$ and $0 \le
\tailr_\nu<\infty$, and for some $\kappa\in[0,2\alpha)$,
%
%e2.11 #&#
\begin{equation}
\label{etail-ind2} \nu(x+I) = O \bigl(\overline{\nu}(x)/x^{1-\kappa
} \bigr),\qquad x
\to\infty.
\end{equation}
Suppose Assumptions~\ref{abasic}\textup{(b)} and \textup{(c)}
hold with $\iof(x, T_0)$ being replaced with $\tilde\iof(x,T_0)$
and $L(x) = x^{\rx/2}$, where $\rx\in(0, 2\alpha- \kappa)$ is a
fixed number. Then the SRT \eqref{erenewal} holds for $\renf$,
where $h>0$ is arbitrary if $\nu$ is nonarithmetic, and is the span
of $\nu$ otherwise.
\end{theorem}

%re3 #&#
\begin{remark*}%\rm\rule{1em}{0em}
\begin{longlist}[(1)]%[itemsep=0ex, leftmargin=0ex,
%itemindent=2.2\parindent, label=\arabic*.]
%
\item[(1)]
Since \eqref{etail-ind} implies $\overline{\cdf}(x)\sim1/\rvf(x)$ and
$\cdf(-x) \sim\tailr_\nu/\rvf(x)$ as $x\to\infty $ (\cite{bingham89},
Theorem~8.2.1), as in Theorem~\ref{trenewal}, once Lemma~\ref{lsmall-n} is established, the rest of the proof of Theorem~\ref{trenewal-ind} is standard.

\item[(2)] Condition \eqref{etail-ind2} can be written as
$x\nu(x+I)/\overline{\nu}(x) = O(x^\kappa)$. Therefore,
it is analogous to \eqref{etail-power} in Proposition~\ref{pprior}. Indeed, our proof of Theorem~\ref{trenewal-ind}
will rely on Proposition~\ref{pprior}.
\end{longlist}
\end{remark*}

%s3 #&#
\section{Proofs for SRT} \label{sproof-renewal}
We shall always denote $M_n = \max_{1\le i\le n} X_i$, and $J =
(-h,h]$. Note $I - I = (-h,h)\subset J$.

%s3.1 #&#
\subsection{An auxiliary result} \label{ssaux}
Some of the notation and arguments in this subsection will also be
used in the proof of Proposition~\ref{pprior}. First, observe that
since for any $x>0$ and $y$, there are at most two $x+k h + J$,
$k\in\Ints$, that contain $y$, then for any $n\ge0$ and event $E$,
%
%e3.1 #&#
\begin{eqnarray}
\label{eoverlap} \sum_k \mathbb{P}\{S_n
\in x + k h + J, E\} &=& \mean \biggl[\sum_k
\mathbf{1}\{S_n \in x + k h + J\} I_E \biggr]
\nonumber
\\[-8pt]
\\[-8pt]
\nonumber
&\le&\mean(2 I_E) = 2\mathbb{P}(E).
\end{eqnarray}

Let $Y_1, Y_2, \ldots $ be i.i.d. following the distribution of $X$
conditional
on $X>0$ and denote
\[
S_n^\pm= \sum_{i=1}^n
X_i^\pm, \qquad N_n = \sum
_{i=1}^n\mathbf{1}\{X_i>0\},\qquad
V_n = \sum_{i=1}^n
Y_i, \qquad\tilde M_n = \max_{i\le n}
Y_i.
\]
Then $S_n = S_n^+ - S_n^-$ and
\[
\mathbb{P}\{Y_i > x\} = \overline{\cdf} \bigl(x^+ \bigr)/p^+ \sim
1/\tilde \rvf(x)\qquad \mbox{with } \tilde\rvf(x) = p^+ \rvf(x).
\]
For $n\ge1$ and $x>0$, define
%
%e3.2 #&#
\begin{equation}
\label{ezeta-gamma} \zeta_{n,x} = a_n^{1-\gamma}
x^\gamma\qquad \mbox{where $(1+\alpha)/(1+2\alpha)<\gamma<1$}
\end{equation}
and
\begin{eqnarray*}
E_{n,x}\sp{(3)} &=& \bigl\{S_n^+\in x+I, X_i>
\zeta_{n,x} \mbox{ for at least two } i=1,\ldots,n \bigr\},
\\
E_{n,x}\sp{(2)} &=& \bigl\{S_n^+\in x+I,
M_n>x/2 \bigr\}\setminus E_{n,x}\sp{(3)},
\\
E_{n,x}\sp{(1)} &= &\bigl\{S_n^+\in x+I,
\zeta_{n,x}< M_n\le x/2 \bigr\}\setminus E_{n,x}
\sp{(3)},
\\
E_{n,x}\sp{(0)} &=& \bigl\{S_n^+\in x+I \bigr\}\,\Big\backslash\,
\bigcup_{i=1}^3 E_{n,x}\sp{(i)}.
\end{eqnarray*}

Let $C_\cdf, C_\cdf', \ldots$ denote constants that only depend on
$\cdf$ (and possibly the fixed~$h$) and may change from line to line.
The auxiliary result we need is the following:

%
%le3.1 #&#
\begin{lemma} \label{lE-bound}
Let \eqref{etail-svf} and \eqref{etail-ratio-2} hold. Fix
$\delta\in(0,1)$ such that $\delta^{\gamma\alpha/2} < p^+$ and
$\delta^{1-\gamma} < 1/2$. Let $p=1$ if $p^+=1$, or $p =9 p^+/10$
if $p^+\in(0,1)$. Then for all $x\gg1$, $n_0\gg1$, $n_0 \le
n\le\rvf(\delta x)$, $p n\le m\le n$ and $T\ge1$,
%
%e3.3 #&#
\begin{equation}
\label{eE3-2-bound} \mathbb{P}\bigl\{E_{n,x}\sp{(i)} | N_n = m
\bigr\} \le \frac{C_\cdf n \overline{\cdf}(x)}{x} \biggl[ T+ \frac{\iof(2x,
T/3)}{a_n} \biggr],\qquad  i=3,2
\end{equation}
and
%
%e3.4 #&#
\begin{equation}
\label{eE1-0-bound} \mathbb{P}\bigl\{E\sp{(i)}_{n,x} | N_n = m
\bigr\} \le \frac{C_\cdf T n \overline{\cdf}(x)}{x},\qquad i=1,0.
\end{equation}
\end{lemma}

\begin{pf}
Notice that for $n\le\rvf(\delta x)$, $\zeta_{n,x} \le(\delta
x)^{1-\gamma} x^\gamma< x/2$ and $n \le\zeta_{n,x}$. Conditional
on $N_n=m$, $S_n^+ \sim V_m$. Then for $1\le m\le n$,
\begin{eqnarray*}
&& \mathbb{P}\bigl\{E_{n,x}\sp{(3)} | N_n = m\bigr\}
\\
&&\qquad\le m^2 \mathbb{P}\{V_m\in x+I, Y_{m-1}>
\zeta_{n,x}, Y_m> \zeta_{n,x}\}
\\
&&\qquad\le n^2 \sum_{k=0}^\infty
\mathbb{P}\{V_m \in x + I, Y_{m-1}>\zeta_{n,x},
Y_m \in \zeta_{n,x} + k h + I\}
\\
&&\qquad\le n^2 \sum_{k=0}^\infty
\mathbb{P}\{V_{m-1}\in x - \zeta_{n,x} - k h + J,
Y_{m-1}> \zeta_{n,x}, Y_m\in
\zeta_{n,x} + k h + I\},
\end{eqnarray*}
where the last line is due to $I-I\subset J$. Then by independence
of $Y_i$, the last inequality yields
%
%e3.5 #&#
\begin{equation}
\label{eevents} \mathbb{P}\bigl\{E_{n,x}\sp{(3)} | N_n=m
\bigr\} \le \frac{n^2}{p^+} \sum_{k=0}^\infty
\mathbb{P}\{X\in\zeta_{n,x} + k h + I\} Q_k,
\end{equation}
where
\[
Q_k = Q_k(m,n,x) = \mathbb{P}\{V_{m-1}\in x -
\zeta_{n,x} - k h + J, Y_{m-1}> \zeta_{n,x}\}.
\]
To bound the RHS of \eqref{eevents}, let
\[
D_k = D_k(n,x,T) = \bigl[ \omega_I(
\zeta_{n,x}+k h) - T \bigr]^+.
\]
Then for $k\ge0$,
\begin{eqnarray*}
\mathbb{P}\{X\in\zeta_{n,x} + k h + I\} &\le&\frac{\overline{\cdf}
(\zeta_{n,x} + k h)}{\zeta_{n,x}+k h}(T +
D_k)
\\
&\le&\frac{\overline{\cdf}(\zeta_{n,x})}{\zeta_{n,x}}(T + D_k).
\end{eqnarray*}
Next, for $x\gg1$, $\zeta_{n,x}\ge\zeta_{1,x}> h$. Then for $k\ge
x/h$, $x-\zeta_{n,x} - k h + h < x - k h<0$, implying $Q_k(x)=0$.
Meanwhile, by \eqref{eoverlap},
\[
\sum_{k=0}^\infty Q_k \le2
\mathbb{P}\{Y_{m-1}> \zeta_{n,x}\} = \frac{2 \overline{\cdf}(\zeta_{n,x})}{p^+}.
\]
Combining \eqref{eevents} and the above bounds,
%
%e3.6 #&#
\begin{eqnarray}
\label{eE3-N} \mathbb{P}\bigl\{E_{n,x}\sp{(3)} | N_n=m
\bigr\} &\le &\frac{n^2\overline{\cdf}(\zeta_{n,x})}{p^+\zeta_{n,x}} \sum_{k=0}^\infty
(T + D_k)Q_k
\nonumber
\\[-8pt]
\\[-8pt]
\nonumber
&\le&\frac{2 T n^2\overline{\cdf}(\zeta_{n,x})^2}{(p^+)^2\zeta
_{n,x}} + \frac{n^2\overline{\cdf}(\zeta_{n,x})}{p^+\zeta_{n,x}} \sum_{0\le k<x/h}
Q_k D_k.
\end{eqnarray}
For each $k$,
%
%e3.7 #&#
\begin{eqnarray}
\label{eQ-bound} Q_k &=& \int_{(\zeta_{n,x}, \infty)} \mathbb{P}
\{V_{m-2}\in x - \zeta_{n,x} - k h - z+ J\} \mathbb{P}\{ X\in
\mathrm{d}z | X>0\}
\nonumber
\\[-8pt]
\\[-8pt]
\nonumber
&\le&\frac{\overline{\cdf}(\zeta_{n,x})}{p^+} \sup_t \mathbb{P}\{
V_{m-2} \in t + J\}.
\end{eqnarray}
By the LLTs (\cite{bingham89}, Theorem~8.4.1--2) and the boundedness of
the density $g$, for all $n\gg1$ and $p n \le m \le n$, $\sup_t
\mathbb{P}\{V_{m-2} \in t + J\} \le C_\cdf/\tilde\rvf(m) \le C_\cdf'/a_n$.
Consequently, by
\eqref{eQ-bound}
\[
\sum_{0\le k<x/h} Q_k D_k \le
\frac{C_\cdf\overline{\cdf}(\zeta_{n,x})}{a_n} \sum_{0\le k<x/h} D_k.
\]
Then by \eqref{eE3-N},
%
%e3.8 #&#
\begin{equation}
\label{eE3-D} \mathbb{P}\bigl\{E_{n,x}\sp{(3)} | N_n=m
\bigr\} \le \frac{2T n^2 \overline{\cdf}(\zeta_{n,x})^2}{(p^+)^2 \zeta_{n,x}} + \frac{C_\cdf n^2 \overline{\cdf}(\zeta_{n,x})^2}{a_n p^+ \zeta_{n,x}} \sum
_{0\le k<x/h} D_k.
\end{equation}

Observe that $D_k = [\omega_I(\zeta_{n,x}+k h) - T]^+\le
[\omega_J(y)-T]^+$ for $y\in\zeta_{n,x} + k h +I$. Then for $x\gg
1$,
\begin{eqnarray*}
\sum_{0\le k<x/h} D_k &\le&\frac{1}{h}
\sum_{0\le k<x/h} \int_{\zeta_{n,x} + k h}^{\zeta_{n,x} + k h+h}
\bigl[\omega_J(y)-T \bigr]^+ \,\mathrm{d}y
\\
&\le&\frac{1}{h} \int_0^{2x} \bigl[
\omega_J(y) - T \bigr]^+ \,\mathrm{d}y.
\end{eqnarray*}
Set $y_0>0$, such that for $y\ge y_0$,
\begin{eqnarray*}
\omega_J(y) &=& \frac{y \mathbb{P}\{X\in y-h+I\}}{\overline{\cdf
}(y)} + \frac{y\mathbb{P}\{X\in y+I\}}{\overline{\cdf}(y)}
\\
&\le&2\omega_I(y-h) + \omega_I(y).
\end{eqnarray*}
On $[0,y_0]$, $\omega_J(y)\le y_0/\overline{\cdf}(y_0)$. By
$[2\omega_I(y-h) + \omega_I(y)-T]^+ \le2[\omega_I(y-h)-T/3]^+ +
[\omega_I(y)-T/3]^+$,
%
%e3.9 #&#
\begin{eqnarray}
\label{eD-bound} \sum_{0\le k<x/h} D_k &\le&
\frac{1}{h}\int_0^{y_0}
\omega_J(y) \,\mathrm{d}y + \frac{1}{h}\int
_{y_0}^{2x} \bigl[\omega_J(y) - T
\bigr]^+ \,\mathrm{d}y
\nonumber
\\[-8pt]
\\[-8pt]
\nonumber
&\le&\frac{C_\cdf'+C_\cdf\iof(2x, T/3)}{h}.
\end{eqnarray}
This combined with \eqref{eE3-D} yields for all $T\ge1$,
\[
\mathbb{P}\bigl\{E_{n,x}\sp{(3)} | N_n=m\bigr\} \le
\frac{C_\cdf n^2 \overline
{\cdf}
(\zeta_{n,x})^2}{\zeta_{n,x}} \biggl[ T+ \frac{\iof(2x, T/3)}{a_n} \biggr].
\]
By \cite{doney97}, page~462, for $x\gg1$ and $n\le\rvf(x)$, $n
\overline{\cdf}(\zeta_{n,x})^2/\zeta_{n,x} \le C_\cdf\overline
{\cdf}(x)/x$.
Insert this inequality into the above one. Then
\eqref{eE3-2-bound} follows for $E_{n,x}\sp{(3)}$.

Since $E_{n,x}\sp{(2)} = \{S_n^+\in x+I$, one $X_i^+>x/2$, all other
$X_i^+\le\zeta_{n,x}\}$, then
\begin{eqnarray*}
&& \mathbb{P}\bigl\{E_{n,x}\sp{(2)} | N_n=m\bigr\}
\\
&&\qquad \le m \mathbb{P}\{V_m \in x + I, \tilde M_{m-1}\le
\zeta_{n,x}, Y_m>x/2\}
\\
&&\qquad\le m \sum_{k=0}^\infty \mathbb{P}
\{V_m \in x + I, \tilde M_{m-1}\le \zeta_{n,x},
Y_m \in x/2 + k h + I\}
\\
&&\qquad\le n \sum_{k=0}^\infty \mathbb{P}
\{V_{m-1}\in x/2 - k h + J, \tilde M_{m-1} \le
\zeta_{n,x}, Y_m \in x/2 + k h + I\}.
\end{eqnarray*}
Denote $Q_k' = Q_k'(m, n,x) = \mathbb{P}\{V_{m-1}\in x/2 - k h + J,
\tilde M_{m-1}\le\zeta_{n,x}\}$ and
\[
D_k' = D_k'(n,x,T) = \bigl[
\omega_I(x/2+k h) - T \bigr]^+.
\]
Then as the argument for \eqref{eE3-N},
\begin{eqnarray*}
\mathbb{P}\bigl\{E_{n,x}\sp{(2)} | N_n=m\bigr\} &\le&
\frac{C_\cdf n \overline
{\cdf}
(x)}{x} \sum_{k=0}^\infty \bigl(T
+ D_k' \bigr)Q_k'
\\
&\le&\frac{C_\cdf n\overline{\cdf}(x)}{x} \biggl( 2 T + \sum_{0\le k<x/h}
Q_k' D_k' \biggr).
\end{eqnarray*}
By the LLTs, $Q_k' \le\sup_t \mathbb{P}\{V_{m-1}\in t + J\} \le
C_\cdf/\tilde\rvf(m) \le C_\cdf'/a_n$. On the other hand,
$\sum_{0\le k < x/h} D_k'$ has the same bound \eqref{eD-bound}.
Then \eqref{eE3-2-bound} follows for $E_{n,x}\sp{(2)}$.

To finish the proof, we need the next general result, which is
essentially due to~\cite{denisov08ap}; see also \cite
{doney97,gouezel11} for results restricted to the arithmetic or operator
cases.

%
%le3.2 #&#
\begin{lemma}[(Denisov, Dieker and Shneer \cite{denisov08ap})]
\label{llocal-ldp}
Let \eqref{etail-svf} and \eqref{etail-ratio-2} hold. There are
$C_\cdf>0$ and $C_\cdf'>0$, such that for any positive sequence
$s_n\to\infty $,
\begin{eqnarray}
\mathbb{P}\{S_n \in x+I, M_n\le s_n\} \le
C_\cdf'(1/s_n + 1/a_n)
e^{-x/s_n + C_\cdf n/\rvf(s_n)},
\nonumber\\
\eqntext{\mbox{all $x>0$ and $n\gg1$}.}
\end{eqnarray}
\end{lemma}

Continuing the proof of Lemma~\ref{lE-bound}, since $E_{n,x}\sp{(1)}
= \{S_n^+\in x+I$, one $X_i^+ \in(\zeta_{n,x},
x/2]$, all other $X_i^+\le\zeta_{n,x}\}$, then
\begin{eqnarray*}
&& \mathbb{P}\bigl\{E_{n,x}\sp{(1)} | N_n = m\bigr\}
\\
&&\qquad\le m \mathbb{P}\{V_m \in x + I, \tilde M_{m-1}\le
\zeta_{n,x} <Y_m \le x/2\}
\\
&&\qquad= m\int_{(\zeta_{n,x}, x/2]} \mathbb{P}\{V_{m-1} \in x-z + I,
\tilde M_{m-1}\le\zeta_{n,x}\} \\
&&\hspace*{48pt}\qquad\quad{}\times\mathbb{P}\{X\in\mathrm{d}z |
X>0\}
\\
&&\qquad\le\frac{n \overline{\cdf}(\zeta_{n,x})}{p^+} \sup_{t\ge x/2} \mathbb {P}
\{V_{m-1}\in t+I, \tilde M_{m-1} \le\zeta_{n,x}\}.
\end{eqnarray*}
Since $p n \le m\le n \le\rvf(\zeta_{n,x}) \le\tilde
\rvf(\zeta_{n, x})$, applying Lemma~\ref{llocal-ldp} to
$\mathbb{P}\{V_{m-1}\in t + I, \tilde M_{m-1}\le\zeta_{n,x}\}$ for
$t\ge
x/2$, with $s_n = \zeta_{n,x}$,
\[
\mathbb{P}\bigl\{E_{n,x}\sp{(1)} | N_n = m\bigr\} \le
C_\cdf n \overline{\cdf }(\zeta_{n,x}) e^{-x/(2\zeta_{n,x})}/a_n,\qquad
p n \le m\le n.
\]
By $n \overline{\cdf}(\zeta_{n,x}) \sim\rvf(a_n)/ \rvf(\zeta_{n,x})
\le
1$ and $e^{-x/(2\zeta_{n,x})}/ a_n \le C_\cdf n\overline{\cdf}(x)/x$
(cf.~\cite{doney97,gouezel11}), \eqref{eE1-0-bound} follows for
$E_{n,x}\sp{(1)}$. Finally, by Lemma~\ref{llocal-ldp}, for $p n
\le m \le n$,
\[
\mathbb{P}\bigl\{E_{n,x}\sp{(0)} | N_n = m\bigr\} =
\mathbb{P}\{V_m \in x + I, \tilde M_m\le
\zeta_{n,x}\} \le C_\cdf e^{-x/\zeta_{n,x}}/a_n,
\]
and \eqref{eE1-0-bound} follows for $E_{n,x}\sp{(0)}$.
\end{pf}

\begin{pf*}{Proof of Lemma~\ref{llocal-ldp}}
This essentially is Lemma~7.1(iv) combined with Proposition~7.1 in
\cite{denisov08ap}. That lemma assumes $s_n$ to be some
specific sequence and $\cdf(-x)$ to be regularly varying at
$\infty$. Both assumptions can be removed. To start with,
for any distribution $\cdf$ and $s>0$ with $\cdf(s)>0$, define
$\tilde\cdf(\mathrm{d}x) = e^{-\psi(s) + x/s} \mathbf{1}\{x\le s\}
\cdf(\mathrm{d}x)$,
where
\[
\psi(s) = \ln\mean \bigl[e^{X/s}\mathbf{1}\{X\le s\} \bigr]\qquad \mbox
{with } X\sim \cdf.
\]
Let $S_n = X_1+\cdots+X_n$ with $X_i$ i.i.d. $\sim\cdf$ and $\tilde S_n
= \tilde X_1+\cdots+\tilde X_n$ with $\tilde X_i$ i.i.d. $\sim\tilde
\cdf$.
Then
\[
\mathbb{P}\{S_n \in x+I, M_n\le s\} \le e^{-x/s+n\psi(s)}
\mathbb {P}\{ \tilde S_n \in x + I\}.
\]
By $\ln\mean Z = \ln[1+\mean(Z-1)] \le\mean(Z-1)$ for any $Z\ge0$
and $e^x-1\le2 x$ for $x\le1$, the following bounds hold:
\begin{eqnarray*}
\psi(s) &\le&\mean \bigl[e^{X/s}\mathbf{1}\{X\le s\} -1 \bigr]
\\
&\le&\mean \bigl[ \bigl(e^{X/s}-1 \bigr)\mathbf{1}\{X\le s\} \bigr]
\\
&\le&\mean \bigl[ \bigl(e^{X/s}-1 \bigr)\mathbf{1}\{0<X\le s\} \bigr]
\\
&\le&2 s^{-1} \mean\bigl[X\mathbf{1}\{0<X\le s\}\bigr].
\end{eqnarray*}
By integration by parts and Karamata's theorem (\cite{bingham89},
Theorem~1.5.11), \eqref{etail-svf} alone implies that for $p\ge1$,
%
%e3.10 #&#
\begin{eqnarray}
\label{eKaramata-p} \int_0^s u^p \cdf(
\mathrm{d}u) &=& p \int_0^s \overline{\cdf}(u)
u^{p-1} \,\mathrm{d}u - \overline{\cdf}(s) s^p
\nonumber
\\[-8pt]
\\[-8pt]
\nonumber
&\sim&\frac{\alpha s^p}{(p-\alpha) \rvf(s)}\to\infty, \qquad s\to\infty
\end{eqnarray}
and hence $\psi(s) \le2 s^{-1} \mean[X\mathbf{1}\{0<X\le s\}] \sim
C_\cdf/\rvf(s)$. Let $\tilde S_n$ be defined with $s=s_n$. Then
$\mathbb{P}\{S_n \in x+I, M_n\le s_n\} \le e^{-x/s_n + C n /\rvf(s_n)}
\mathbb{P}\{\tilde S_n \in x + I\}$, and so it only remains to check
%
%e3.11 #&#
\begin{equation}
\label{etilted-S} \mathbb{P}\{\tilde S_n\in x + I\} \le
C_\cdf(1/s_n + 1/a_n).
\end{equation}
Since \eqref{etail-ratio-2} holds as well, there is $s_0>0$, such
that for $s>s_0$,
%
%e3.12 #&#
\begin{eqnarray}
\label{eKaramata-n} \int_{-s}^0 |u|^p
\cdf(\mathrm{d}u) &\le &p \int_0^s \cdf(-u)
u^{p-1} \,\mathrm{d}u
\nonumber
\\
&\le& C_\cdf+ \tailr p \int_{s_0}^s
\overline{\cdf}(u) u^{p-1} \,\mathrm{d}u
\\
&\le& C_\cdf s^p/\rvf(s).\nonumber
\end{eqnarray}
Let $\mu_p(s):= \mean[|X|^p\mathbf{1}\{|X|\le s\}]$. Then for $s\gg1$,
by \eqref{eKaramata-p} and \eqref{eKaramata-n},
%
%e3.13 #&#
\begin{equation}
\label{emu-p} C_\cdf s^p / \rvf(s) \le\mean
\bigl[X^p \mathbf{1}\{0<X\le s\} \bigr] \le\mu_p(s) \le
C_\cdf' s^p/\rvf(s).
\end{equation}
It follows that $\limsup_{x\to\infty } x^2 \overline{G}(x)/\mu
_2(x)<\infty$, where $G(x) = \mathbb{P}\{|X|\le x\}$, so by
Proposition~7.1 in \cite{denisov08ap}, for all $n\gg1$,
%
%e3.14 #&#
\begin{equation}
\label{etilted-S-2} \sup_x \mathbb{P}\{\tilde S_n
\in x + I\} \le C_\cdf(1/s_n + 1/r_n),
\end{equation}
where $r_n>0$ is the solution to $Q(x):= x^{-2} \mu_2(x) + \overline
{G}(x) = 1/n$, which exists and is unique for all $n\gg1$. On the
one hand, since $Q(x)\ge\overline{\cdf}(x) \sim1/\rvf(x)$, $r_n
\ge
C_\cdf a_n$. On the other, by \eqref{emu-p}, $Q(x)\le C_\cdf
/\rvf(x)$ and then $r_n\le C_\cdf' a_n$. Then \eqref{etilted-S}
follows from \eqref{etilted-S-2}.
\end{pf*}

%s3.2 #&#
\subsection{Proof of Theorem \texorpdfstring{\protect\ref{trenewal}}{2.1}} \label{sssmall-n}
We need two lemmas for the proof of Theorem~\ref{trenewal}.

%le3.3 #&#
\begin{lemma} \label{lsmall-n}
Let \eqref{etail-svf} and \eqref{etail-ratio-2} hold. Then
Assumption~\ref{abasic} implies
\[
\lim_{\delta\to0+} \limsup_{x\to\infty }
\frac{x}{\rvf(x)} \sum_{n\le\rvf(\delta x)} \mathbb{P}
\{S_n\in x + I\} = 0.
\]
\end{lemma}

%le3.4 #&#
\begin{lemma} \label{lmiddle-n}
Let \eqref{etail-svf} and \eqref{etail-ratio} hold. Given
$0<\delta<1$, let $J_\delta(x) = (\rvf(\delta x), \rvf(x/\delta))$.
Then
%
%e3.15 #&#
\begin{equation}
\label{esrt-truncate} \lim_{x\to\infty } \frac{x}{\rvf(x)}\sum
_{n\in J_\delta(x)} \mathbb{P}\{S_n \in x+I\} = \alpha h\int
_\delta^{1/\delta} x^{-\alpha} g(x) \,\mathrm{d}x.
\end{equation}
\end{lemma}

Assume the lemmas are true for now.
Since $\mathbb{P}\{S_n\in x + I\} =
O(1/a_n)$ and
\[
\sum_{n\ge\rvf(x/\delta)} 1/a_n \sim
\frac{\rvf(x/\delta)}{(\alpha^{-1} - 1) x/\delta} \sim\delta ^{1-\alpha} \frac{\rvf(x)}{(\alpha^{-1} - 1) x},
\]
by \eqref{esrt-truncate},
\[
\lim_{x\to\infty } \frac{x}{\rvf(x)}\sum
_{n>\rvf(\delta x)} \mathbb{P}\{S_n \in x+I\} = \alpha h\int
_0^{1/\delta} x^{-\alpha} g(x) \,\mathrm{d}x.
\]
Combining this with Lemma~\ref{lsmall-n} and letting $\delta\to
0+$, we then get \eqref{erenewal}.

\begin{pf*}{Proof of Lemma~\ref{lsmall-n}}
Denote $\Omega_{n,x} = \{S_n \in x + I\}$. By Assumption~\ref{abasic}, it suffices to show
%
%e3.16 #&#
\begin{equation}
\label{elow-cut-L} \lim_{\delta\to0+} \limsup_{x\to\infty }
\frac{x}{\rvf(x)} \sum_{L(x) \le n\le\rvf(\delta x)} \mathbb{P}(
\Omega_{n,x}) = 0.
\end{equation}
Let $\delta>0$ such that $\delta^{-\gamma\alpha/2} p^+>1 > 2
\delta^{1-\gamma}$. Set $p =1$ if $p^+=1$, and $p=9 p^+/10$ if
$p^+<1$. Then
\[
\mathbb{P}(\Omega_{n,x}) \le\sum_{p n\le m\le n}
\mathbb{P}( \Omega_{n,x} | N_n = m) \mathbb{P}
\{N_n=m\} + \mathbb{P}\{N_n < p n\}.
\]
For each $p n\le m\le n$, conditional on $N_n=m$, $S_n^+$ and
$S_n^-$ are independent. Therefore,
\begin{eqnarray*}
& &\mathbb{P}(\Omega_{n,x} | N_n = m)
\\
&&\qquad= \int_0^\infty \mathbb{P}\bigl
\{S_n^+ \in x+z+I | N_n = m\bigr\} \mathbb{P}\bigl
\{S_n^- \in\mathrm{d}z | N_n = m\bigr\}.
\end{eqnarray*}
Since $\{S_n^+ \in x+z+I\} = \bigcup_{i=0}^4 E_{n,x+z}\sp{(i)}$, by Lemma~\ref{lE-bound}, for $x\gg1$, $L(x)\le n\le\rvf(\delta x)$
and $p n \le m\le n$,
\begin{eqnarray*}
& &\mathbb{P}(\Omega_{n,x} | N_n = m)
\\
&&\qquad\le C_\cdf n\int_0^\infty
\frac{\overline{\cdf}(x+z)}{x+z} \biggl[T_0+ \frac{\iof(2x+2z, T_0)}{a_n} \biggr] \mathbb{P}
\bigl\{S_n^-\in\mathrm {d}z | N_n=m\bigr\}
\\
&&\qquad= C_\cdf n \mean \biggl\{ \frac{\overline{\cdf}(x_n)}{x_n} \biggl[T_0
+ \frac{\iof(2x_n, T_0)}{a_n} \biggr] \Big| N_n=m \biggr\},
\end{eqnarray*}
where $x_n = x + S_n^-$. $N_n$ is the sum of $n$ independent
Bernoulli random variables each with mean $p^+$. If $p^+\in(0,1)$,
then by Chernoff's inequality, for $n\gg1$, $\mathbb{P}\{N_n < p n\}
\le
e^{-\lambda n}$, where $\lambda=\lambda(p^+)>0$ is a constant; cf. \cite{taovu06}, Corollary~1.9. As a result,
%
%e3.17 #&#
\begin{equation}
\label{eE-n} \mathbb{P}(\Omega_{n,x}) \le C_\cdf n \mean
\biggl\{ \frac{\overline{\cdf}(x_n)}{x_n} \biggl[T_0 +\frac{\iof(2x_n,
T_0)}{a_n} \biggr]
\biggr\} + e^{-\lambda n}.
\end{equation}
If $p^+=1$, then $N_n\equiv n$ and $S_n^-=0$, so by setting
$\lambda=\infty$, the above inequality still holds.

Since $x_n\ge x$, given $c>1$, $\overline{\cdf}(x_n)/x_n\le
\overline{\cdf}(x)/x\le c/[x \rvf(x)]$ for $x\gg1$. Then, writing
\[
R(x) = \frac{x}{\rvf(x)} \sum_{L(x)\le n\le\rvf(\delta x)}
\frac{n}{a_n} \mean \biggl[ \frac{\iof(2x_n, T_0)}{x_n\rvf(x_n)} \biggr],
\]
we have
\begin{eqnarray*}
&& \frac{x}{2\rvf(x)}\sum_{L(x)\le n\le\rvf(\delta x)} n \mean \biggl\{
\frac{\overline{\cdf}(x_n)}{x_n} \biggl[T_0 +\frac{\iof(2x_n,
T_0)}{a_n} \biggr] \biggr\}
\\
&&\qquad\le\frac{T_0}{\rvf(x)^2}\sum_{L(x)\le n\le\rvf(\delta x)} n + R(x).
\end{eqnarray*}
The first term on the RHS is $O(\rvf(\delta x)^2/\rvf(x)^2) =
\delta^{2\alpha}$. To bound $R(x)$, write $\theta= 1/\alpha-1$.
Then $\theta>0$ and $n/a_n = n^{-\theta}/\svf(n)$. Consider two
cases.

\textit{Case} 1: $\alpha\in(0,1/2)$. Then $\theta>1$. By Assumption~\ref{abasic},
\[
\iof(2 x_n, T_0) = o \bigl(\rvf(2 x_n)^2/u_\theta(2
x_n) \bigr),
\]
and since $L\in\CR_c$ with $c \in[0,\alpha]$, $u_\theta\in
\CR_{c(1-\theta)}$. As a result,
\begin{eqnarray*}
R(x) &\le&\frac{x}{\rvf(x)} \biggl(\sum_{n\ge L(x)}
\frac{n}{a_n} \biggr) \max_{n\ge L(x)} \mean \biggl[
\frac{\iof(2x_n,
T_0)}{x_n \rvf(x_n)} \biggr]
\\
&=& o \biggl( \frac{x u_\theta(x)}{\rvf(x)} \max_{n\ge L(x)} \mean \biggl[
\frac{\rvf(x_n)}{x_n u_\theta(x_n)} \biggr] \biggr), \qquad x\to\infty.
\end{eqnarray*}
Since $\rvf(x)/x u_\theta(x)\in\CR_b$ with
\[
b = \alpha-1 + c(\theta-1) \le\alpha-1 + \alpha(1/\alpha-2) = -\alpha<0,
\]
then $\mean[\rvf(x_n)/x_n u_\theta(x_n)] = O(\rvf(x)/x
u_\theta(x))$. It follows that $R(x)\to0$ as \mbox{$x\to\infty $}.

\textit{Case} 2: $\alpha=1/2$. Since $x/\rvf(x) \le2 x_n/\rvf(x_n)$
for $x\gg1$,
\[
R(x) = O(1) \sum_{L(x)\le n\le\rvf(\delta x)} \frac{n}{a_n} \mean
\biggl[\frac{\iof(2 x_n, T_0)}{\rvf(x_n)^2} \biggr].
\]
If $\tilde u(x)/\tilde u(L(x))\to1$, then by Assumption~\ref{abasic},
\[
\iof(2 x_n, T_0)/\rvf(x_n)^2 =
O \bigl(1/\tilde u(2 x_n) \bigr).
\]
Since $\tilde u(x)$ is increasing, then
\begin{eqnarray*}
R(x) &=& O \bigl(1/\tilde u(x) \bigr) \sum_{\rvf(L(x)) \le n \le\rvf
(\delta x)}
\frac{n^{-1}}{\svf(n)}
\\
&=& O \bigl(1/\tilde u(x) \bigr) \int_{\rvf(L(x))}^{\rvf(x)}
\frac{y^{-1}}{\svf(y)} \,\mathrm{d}y = o(1).
\end{eqnarray*}
If $\tilde u(x)/\tilde u(L(x))\not\to1$, then by Assumption~\ref{abasic},
\[
\iof(2 x_n, T_0)/\rvf(x_n)^2 =
o \bigl(1/\tilde u(2 x_n) \bigr) = o \bigl(1/\tilde u(x) \bigr),
\]
and hence
\[
R(x) = o \bigl(1/\tilde u(x) \bigr) \sum_{n \le\rvf(\delta x)}
\frac{n^{-1}}{\svf(n)} = o(1).
\]

Thus, for all $\alpha\in(0,1/2]$, $R(x)= o(1)$. Finally, if
$p^+\in(0,1)$, then given $c>0$ such that $\lambda c> 1-\alpha$,
for $x\gg1$, $\sum_{n\ge L(x)} e^{-\lambda n} \le\sum_{n\ge c\ln
x} e^{-\lambda n}= O(x^{-\lambda c}) = o(\rvf(x)/x)$. Then by
summing \eqref{eE-n} over $L(x)\le n \le\rvf(\delta x)$ and taking
the limit as $x\to\infty $ followed by $\delta\to0$, the proof is
complete.
\end{pf*}

\begin{pf*}{Proof of Lemma~\ref{lmiddle-n}}
If $X$ is arithmetic or nonlattice, then \eqref{esrt-truncate} is
well-known \cite
{garsialamperti63,doney97,erickson70tams,vatutin13tpa}. The only
remaining case is where $X$ is lattice
but nonarithmetic. While Theorem~2 in \cite{erickson70tams} correctly
states that \eqref{esrt-truncate} still holds in this case, the
argument therein cannot establish the fact as it overlooks issues
caused by the discrete nature of $X$.

Let $X$ be concentrated in $a + d\Ints$ with $a/d>0$ being
irrational and $d>0$ the span. If $h\ge d$, then choose $k\in
\Nats$ such that $h'=h/k<d$. Letting $I' = (0,h']$ and $x_j = x+j
h'$, $\mathbb{P}\{S_n\in x+I\} = \sum_{j=0}^{k-1} \mathbb{P}\{S_n\in
x_j + I'\}$,
$x_j/\rvf(x_j) \sim x/\rvf(x)$, and hence if \eqref{esrt-truncate}
holds for $\mathbb{P}\{S_n\in x+I'\}$, it holds for $\mathbb{P}\{S_n
\in x + I\}$ as
well. Thus, without loss of generality, let $0< h <d$.

For $z\in\Reals$, denote $L_z:=d\Ints\cap(z + I)$. Since
$h<d$, $L_z$ contains at most one point. For $x>0$ and $n\ge1$, if
$L_{x-n a}=\{d k\}$, by Gnedenko's LLT, $\mathbb{P}\{S_n \in x + I\} =
\mathbb{P}\{S_n = n a + d k\} = (d/a_n) [g((n a + k d)/a_n) + o(1)]$ as
$n\to\infty $, where $o(1)$ is uniform in $x$ (\cite{bingham89},
Theorem~8.4.1). Since $g$ has bounded derivative and $|x - (n a + k
d)|<h$, it is seen
\[
\mathbb{P}\{S_n \in x + I\} = \mathbf{1}\{L_{x-n a}\ne
\varnothing\} (d/a_n) \bigl[g(x/a_n) + o(1) \bigr].
\]

For $x\gg1$ and $n\in J_\delta(x)$, $x/a_n \in(\delta, 1/\delta)$.
Since $g>0$ on $[\delta, 1/\delta]$, the above display can be
written as
%
%e3.18 #&#
\begin{eqnarray}
\label{eLLT} %
\mathbb{P}\{S_n \in x + I\} = \mathbf{1}
\{L_{x-n a}\ne\varnothing\} (d/a_n) \bigl[1+
\rx_n(x) \bigr] g(x/a_n),
\nonumber
\\[-8pt]
\\[-8pt]
\eqntext{\mbox{with } \displaystyle\lim_{x\to\infty }\sup_{n\in J_\delta(x)} \bigl|
\rx_n(x)\bigr| = 0.}
\end{eqnarray}

Let $m=m(x)$ and $M = M(x)$ be the smallest and largest integers
in $(\rvf(\delta x),   \rvf(x/\delta) + 1)$, respectively. Fix
integers $m = N_1 < N_2 < \cdots< N_s < N_{s+1} = M$, where $N_i =
N_i(x)$ and $s=s(x)$, such that as $x\to\infty $,
\[
\min_{1\le i\le s} [N_{i+1} - N_i]\to\infty,\qquad
\max_{1\le i\le s} [a_{N_{i+1}} - a_{N_i}] = o(x).
\]
Then by \eqref{eLLT}
\[
\sum_{n\in J_\delta(x)} \mathbb{P}\{S_n \in x+I\}
\sim d \sum_{j=1}^s \sum
_{n=N_j}^{N_{j+1}-1} \mathbf{1}\{L_{x-n a}\ne
\varnothing\} \frac{g(x/a_n)}{a_n}.
\]
By the choice of $N_1, \ldots, N_{s+1}$, for each $1\le j\le s$ and $n=N_j,
\ldots, N_{j+1}-1$,
%
%e3.19 #&#
\begin{eqnarray}
\label{eg-approx} %
g(x/a_n)/a_n = \bigl[1+
\rx_n(x) \bigr] g(x/a_{N_j})/a_{N_j},
\nonumber
\\[-8pt]
\\[-8pt]
\eqntext{\mbox{with }  \displaystyle\lim_{x\to\infty }
\sup_{n\in J_\delta(x)} \bigl|
\rx_n(x)\bigr|=0.}
\end{eqnarray}
Thus
\[
\sum_{n\in J_\delta(x)} \mathbb{P}\{S_n \in x+I\}
\sim d \sum_{j=1}^s \frac{g(x/a_{N_j})}{a_{N_j}}
\sum_{n=N_j}^{N_{j+1}-1} \mathbf{1}
\{L_{x-n a}\ne \varnothing\}.
\]

Denote $K=\{\omega\in\Coms\dvtx|\omega|=1\}$. Then $L_z\ne
\varnothing$
$\iff$ $e^{2\pi\iunit z/d}$ falls into the arc $\Gamma= \{\omega=
e^{2\pi\iunit\theta}\dvtx\theta\in[-h/d,0)\}\subset K$. Let $c=a/d$
and define $T\dvtx K\to K$ as $T(\omega) = \omega e^{-2\pi\iunit c}$.
Let $\omega_j = e^{2\pi\iunit(x - N_j a)/d}$. Then
\[
\sum_{n=N_j}^{N_{j+1}-1} \mathbf{1}
\{L_{x-n a}\ne \varnothing\} = \sum_{n=0}^{N_{j+1} - N_j-1}
\mathbf{1}\bigl\{T^n(\omega_j)\in\Gamma\bigr\}.
\]

Since $c$ is irrational, $T$ is a homeomorphism of $K$ with no
periodic points, that is, for any $\omega\in K$ and $n\in\Nats$,
$T^n(\omega) \ne\omega$. Then by ergodic theory (\cite{walters82},
Section~6.5), for any $f\in C(K)$, $(1/N) \sum_{n=0}^{N-1} f(T^n
\omega)\to\int f \,\mathrm{d}\mu$ uniformly in $\omega\in K$, with
$\mu$
the uniform probability measure on $K$. Since $\mu(\Gamma) = h/d$,
and for any $\rx>0$, there are $f$, $g\in C(K)$ with $0\le f(\omega)
\le
\mathbf{1}\{\omega\in\Gamma\} \le g(\omega)\le1$ such that $0\le
\int(g-f) \,\mathrm{d}\mu< \rx$, then
\begin{eqnarray}
\sum_{n=0}^{N_{j+1} - N_j-1} \mathbf{1}\bigl
\{T^n(\omega_j)\in\Gamma\bigr\} = (N_{j+1} -
N_j) \bigl[1+\rx_j(x) \bigr] (h/d),\nonumber
\\
\eqntext{\mbox{with } \displaystyle \lim_{x\to\infty } \sup_{1\le j\le s} \bigl|
\rx_j(x)\bigr|=0.}
\end{eqnarray}
This combined with the previous two displays and then with
\eqref{eg-approx} yields
\begin{eqnarray*}
\sum_{n\in J_\delta(x)} \mathbb{P}\{S_n \in x+I\} &
\sim& d \sum_{j=1}^s \frac{g(x/a_{N_j})}{a_{N_j}}
(N_{j+1} - N_j) (h/d)
\\
& \sim& h \sum_{j=1}^s \sum
_{n=N_j}^{N_{j+1}-1} \frac{g(x/a_n)}{a_n} = h\sum
_{n\in J_\delta(x)} \frac{g(x/a_n)}{a_n}.
\end{eqnarray*}
Multiply both sides by $x/\rvf(x)$ and let $x\to\infty $. Standard
derivation such as the one on page~366 in \cite{bingham89} then
yields \eqref{esrt-truncate}.
\end{pf*}

%s3.3 #&#
\subsection{Proof of Proposition \texorpdfstring{\protect\ref{pprior}}{2.3}}
Given $0<\rx<2\alpha-\kappa$, fix $c\in(\rx/(2\alpha),1)$ such that
$1+\rx< c(\alpha+1-\kappa) + \alpha$. Since $X_i$ are i.i.d.,
\begin{eqnarray*}
&& \mathbb{P}\bigl\{S_n\in x+I, M_n>x^c\bigr
\}
\\
&&\qquad\le n \mathbb{P}\bigl\{S_n\in x+I, X_n >
x^c\bigr\}
\\
&&\qquad= n \sum_{k=0}^\infty \mathbb{P}\bigl
\{S_n \in x+I, X_n \in x^c+k h + I\bigr\}
\\
&&\qquad\le n \sum_{k=0}^\infty \mathbb{P}\bigl
\{S_{n-1}\in x-x^c-k h + J, X_n \in
x^c + k h+I\bigr\}
\\
&&\qquad= n \sum_{k=0}^\infty \mathbb{P}\bigl
\{S_{n-1}\in x-x^c-k h + J\bigr\} \mathbb{P}\bigl\{X\in
x^c + k h+I\bigr\}.
\end{eqnarray*}
Then by \eqref{eoverlap}, $\mathbb{P}\{S_n\in x+I, M_n>x^c\}\le2n
\sup_{t\ge x^c} \mathbb{P}\{X\in t+I\}$. By assumption, for all $t\ge x^c$,
$\mathbb{P}\{X\in t+I\} \le C/[t^{1-\kappa} \rvf(t)] \le
C/[x^{c(1-\kappa)}
\rvf(x^c)]$, where $C>0$ is a constant that may change from line to
line. Then by the choice of $c$,
%
%e3.20 #&#
\begin{eqnarray}
\label{eprior-1} && x\overline{\cdf}(x) \sum_{n\le x^{\rx/2}}
\mathbb{P}\bigl\{S_n\in x+I, M_n>x^c\bigr\}
\nonumber
\\
&&\qquad\le\frac{C x}{\rvf(x)} \sup_{t\ge x^c} \mathbb{P}\{X\in t+I\} \sum
_{n\le x^{\rx/2}} n
\\
\nonumber& &\qquad\le\frac{C x^{1+\rx}}{x^{c(1-\kappa)} \rvf(x^c) \rvf(x)}
= o(1),\qquad x\to\infty.
\end{eqnarray}
Note that if $S_n\in x+I$ and $M_n\le x^c$, then $n\ge x^{1-c}$. By
$x^{\rx/2} = o(\rvf(x^c))$ and Lemma~\ref{llocal-ldp},
\begin{eqnarray*}
\sum_{n\le x^{\rx/2}} \mathbb{P}\bigl\{S_n\in x+I,
M_n\le x^c\bigr\} &=& \sum_{x^{1-c}\le n\le x^{\rx/2}}
\mathbb{P}\bigl\{S_n\in x+I, M_n\le x^c\bigr
\}
\\
&\le& C \sum_{x^{1-c}\le n\le x^{\rx/2}} \bigl(1/x^c +
1/a_n \bigr) e^{-x^{1-c}}
\\
&\le& o \bigl(e^{-x^{1-c}} \bigr).
\end{eqnarray*}
Then $x\overline{\cdf}(x)\sum_{n\le x^{\rx/2}} \mathbb{P}\{S_n\in x+I,
M_n\le x^c\}
= o(1)$, which together with \eqref{eprior-1} completes the proof.

%s3.4 #&#
\subsection{Proof of Theorem \texorpdfstring{\protect\ref{tladder}}{2.4}}
Since $\alpha\in(0,1)$, it is known that $\rvf^+\in\CR_{\alpha
\varrho}$, that is, $\rvf^+(x)$ is regularly varying with exponent
$\alpha\varrho$ \cite{rogozin71}. Let $\omega_I^+(x)$ and
$\iof^+(x,T)$ denote the functions defined by \eqref{eprob-ratio} and
\eqref{eoverflow} for $\ldh$. By assumption, $\iof(x,T) = O(x^{2
c\alpha})$ for some $c\in[0, \varrho)$ and $T>0$. We shall show
that for any $\gamma\in(c,\varrho)$,
%
%e3.21 #&#
\begin{equation}
\label{eladder-overflow} \iof^+(x,T) = O \bigl(x^{2\gamma\alpha
} \bigr).
\end{equation}
Once this is proved, then the proof follows from Corollary~\ref
{crenewal}. For $t>0$,
\[
\mathbb{P}\{\ldh\in t+I\} = \int_0^\infty
\mathbb{P}\{X\in t+y+I\} \renf ^-(\mathrm{d}y),
\]
where $\renf^-(\mathrm{d}t) = \sum_{n=0}^\infty \mathbb{P}\{\ldh
_n^-\in-\mathrm{d}t\}$
concentrates on $[0,\infty)$, with $\ldh_n^-$ the weak decreasing
latter process of $S_n$ (\cite{feller71}, page~399). Then
\begin{eqnarray*}
\bigl[\omega_I^+(t)-T \bigr]^+ &=& \frac{1}{\mathbb{P}\{H>t\}} \bigl[t
\mathbb{P}\{H\in t+I\} - \mathbb{P}\{H>t\} T\bigr]^+
\\
&\le&\frac{1}{\mathbb{P}\{H>t\}} \int_0^\infty \bigl[t
\mathbb {P}\{X\in t+y+I\} - \overline{\cdf}(t+y) T \bigr]^+ \renf^-(\mathrm{d}y)
\\
&\le&\frac{1}{\mathbb{P}\{H>t\}} \int_0^\infty
\frac{t \overline
{\cdf}(t+y)}{t+y} \bigl[\omega_I(t+y)-T \bigr]^+ \renf^-(
\,\mathrm{d}y).
\end{eqnarray*}
Denote $g_y(t) = t \overline{\cdf}(t+y)/(t+y)$. Then
%
%e3.22 #&#
\begin{eqnarray}
\label{eiof-ladder} %
\iof^+(x, T) \le\sum_{i=1}^4
I_i,
\nonumber
\\[-8pt]
\\[-8pt]
\eqntext{\mbox{with }  \displaystyle I_i = \int_{A_i}
\frac{g_y(t)[\omega_I(t+y)-T]^+}{\mathbb{P}\{\ldh>t\}}
 \,\mathrm{d}t \,\renf ^-(\mathrm{d}y),}
\end{eqnarray}
where $A_1 = \{0\le t\le x < y\}$, $A_2 = \{0\le t < y \le x\}$,
$A_3 = \{M\le y\le t\le x\}$, and $A_4 = \{0<y<M, y\le t\le x\}$,
where $M\gg1$ is a fixed number. Fix $0<\beta<\alpha$.

First, let $\varrho\in(0,1)$. For $(t,y)\in A_1$, $\mathbb{P}\{\ldh
>t\} \ge
\mathbb{P}\{\ldh>x\}$. Let $x\gg1$. Then $g_y(t)\le h_y(t):=t/(t+y)^{1+\beta}$ and
\[
I_1 \le\frac{1}{\mathbb{P}\{\ldh>x\}} \int_0^x
\,\mathrm{d}t \int_x^\infty h_y(t)
\bigl[ \omega_I(t+y) - T \bigr]^+ \renf^-(\mathrm{d}y).
\]
Since for each $y\ge x$, $h_y(t)$ is increasing on $[0,x]$,
\begin{eqnarray*}
I_1 &\le&\frac{1}{\mathbb{P}\{\ldh>x\}} \int_0^x
\,\mathrm{d}t \int_x^\infty h_y(x)
\bigl[ \omega_I(t+y) - T \bigr]^+ \renf^-(\mathrm{d}y)
\\
&\le&\frac{1}{\mathbb{P}\{\ldh>x\}} \int_x^\infty
h_y(x) \iof(x+y, T) \renf^-(\mathrm{d}y)
\\
&\le&\frac{C}{\mathbb{P}\{\ldh>x\}} \int_x^\infty
\frac{x}{(x+y)^{1+\beta- 2 c\alpha}} \renf^-(\mathrm{d}y),
\end{eqnarray*}
where $C>0$ is a constant. Since $\ldh^-$ is in the domain of
attraction of stable law with exponent $\alpha(1-\varrho)$
\cite{doney93ptrf}, $\renf^-(x \,\mathrm{d}u)/x^{\alpha(1-\varrho)}$
converges vaguely to $C u^{\alpha(1-\varrho) - 1}\mathbf{1}\{u>0\}
\,\mathrm{d}u$ as
$x\to\infty $, where $C>0$ is a constant; see \cite{bingham89},
pages~361--363. Therefore, by variable substitute $y = x u$,
\[
I_1 \le\frac{C x^{-\beta+2 c\alpha+\alpha(1-\varrho)}}{
\mathbb{P}\{\ldh>x\}} \int_1^\infty
\frac{u^{\alpha(1-\varrho)-1} \,\mathrm{d}u}{(1+u)^{1+\beta- 2
c\alpha}}.
\]
As long as $0<\alpha-\beta\ll1$, the integral is finite. (Recall
that $c\alpha< \alpha\varrho\le1/2$.) Then, for any $\gamma> c$,
$I_1=O(x^{2 \gamma\alpha})$.

To bound $I_2$, observe $\mathbb{P}\{\ldh>t\} \ge\mathbb{P}\{\ldh
>y\}$ for
$(t,y)\in A_2$. Then by $g_y(t)\le\overline{\cdf}(y)$,
\begin{eqnarray*}
I_2 &\le&\int_0^x
\frac{\renf^-(\mathrm{d}y)}{\mathbb{P}\{\ldh>y\}} \int_0^y
g_y(t) \bigl[\omega_I(t+y)-T \bigr]^+ \,\mathrm{d}t
\\
&\le&\int_0^x \frac{\overline{\cdf}(y) \renf^-(\mathrm
{d}y)}{\mathbb{P}\{\ldh
>y\}} \int
_0^y \bigl[\omega_I(t+y)-T
\bigr]^+ \,\mathrm{d}t
\\
&\le&\int_0^x \frac{\overline{\cdf}(y) \iof(2y, T)}{\mathbb{P}\{
\ldh
>y\}
} \renf^-(
\mathrm{d} y).
\end{eqnarray*}
By assumption, $\overline{\cdf}(y) \iof(2y, T)/\mathbb{P}\{\ldh
>y\} =
O(y^{-\beta+2c\alpha+ \alpha\varrho})$ for any $\beta<\alpha$. Since
$\renf^-((0,x])$ is regularly varying with exponent
$\alpha(1-\varrho)$, the integral is of order $O(x^{2\gamma\alpha})$
for any $\gamma>c$.

Let $M\gg1$ such that $\overline{\cdf}(t+y)/\mathbb{P}\{\ldh> t\} <
k_y(t):=t^{\alpha\varrho}/(t+y)^\beta$ for $(t,y)\in A_3$. If
$\beta\in(\alpha\varrho, \alpha)$, then $k_y(t)$ has maximum value
$C/y^{\beta- \alpha\varrho}$, where $C=C(\beta)>0$ is a constant.
Then
\begin{eqnarray*}
I_3 &\le&\int_{A_3} \frac{C t}{(t+y)y^{\beta-\alpha\varrho}} \bigl[
\omega_I(t+y) - T \bigr]^+ \,\mathrm{d}t \,\renf^-(\mathrm{d}y)
\\
&\le&\int_M^x \frac{C\renf^-(\mathrm{d}y)}{y^{\beta- \alpha
\varrho}} \int
_y^x \bigl[\omega_I(t+y)-T
\bigr]^+ \,\mathrm{d}t
\\
&\le&\iof(2x,T) \int_M^x \frac{C\renf^-(\mathrm{d}y)}{y^{\beta-
\alpha
\varrho}}.
\end{eqnarray*}
The integral is of order $O(x^{\alpha-\beta})$. Then by the
assumption on $\iof$, $I_3 = O(x^{2\gamma\alpha})$ for any
$\gamma>c$.

For $I_4$, since $g_y(t)/\mathbb{P}\{\ldh>t\} \le\overline{\cdf}
(t)/\mathbb{P}\{\ldh>t\}$
is bounded,
\begin{eqnarray*}
I_4 &\le&\int_0^M \renf^-(
\mathrm{d} y) \int_y^x \bigl[
\omega_I(t+y)-T \bigr]^+ \,\mathrm{d}t
\\
&\le&\iof(2x, T) \renf^- \bigl([0,M )\bigr) = O \bigl(x^{2 c\alpha
} \bigr).
\end{eqnarray*}
Combining the above bounds for $I_i$ and \eqref{eiof-ladder}, then
\eqref{eladder-overflow} follows when $\varrho\in(0,1)$.

If $\varrho=1$ and $S_n\to\infty $ a.s., then $\renf^-$ is a finite measure
(\cite{feller71}, pages~395--396). It is then not hard to see the above
bounds for $I_i$ still hold. The proof is then complete.

%s3.5 #&#
\subsection{Proof of Theorem \texorpdfstring{\protect\ref{trenewal-ind}}{2.5}}
From the proof of Theorem~\ref{trenewal}, it suffices to prove
Lemma~\ref{lsmall-n} under the assumptions on the $\levy$ measure
$\nu$ of $X$. We will use several times the fact that $X\sim
Y_1+\cdots+Y_N
+ W$, where $Y_i$, $N$ and $W$ are independent,
\[
Y_i\sim Y\sim G(x) = \mathbf{1}\{x>1\}\nu \bigl((1,x) \bigr)/
\nu_0,
\]
with $\nu_0 = \overline{\nu}(1)$, $N\sim\dpois(\nu_0)$, and for
$\theta
\in
\Reals$,
\begin{eqnarray*}
 \mean \bigl[e^{\iunit\theta W} \bigr]
= \exp \biggl\{ \iunit\mu\theta- \sigma^2\theta^2/2 +
\int \bigl(e^{\iunit\theta u} - 1 - \iunit\theta u\mathbf{1}\bigl\{|u|\le1\bigr\} \bigr)
\mathbf{1}\{u\le1\} \nu(\mathrm{d} u) \biggr\}.
\end{eqnarray*}
Then $\mean[e^{t W}]<\infty$ for any $t>0$ (\cite{sato99},
Theorem~25.17). Write $\zeta_N = Y_1+\cdots+Y_N$, and when $N$ is random,
always assume that it is independent of $Y_i$.

%le3.5 #&#
\begin{lemma} \label{lind}
Let \eqref{etail-ind} and \eqref{etail-ind2} hold. Then
given $c \in(0, 2\alpha- \kappa)$, \eqref{elow-cut}~holds
with $L(x) = x^{c/2}$.
\end{lemma}

By this lemma, it suffices to show
\[
\lim_{\delta\to0+} \limsup_{x\to\infty }
\frac{x}{\rvf(x)} \sum_{L(x)\le n\le\rvf(\delta x)} \mathbb{P}(
\Omega_{n,x}) = 0,
\]
where $L(x)=x^{\rx/2}$ and $\Omega_{n,x}=\{S_n\in x+I\}$. For $n\ge
1$,
\[
S_n \sim\zeta_{N_n} + V_n \qquad\mbox{with }
V_n = W_1+\cdots+W_n,
\]
where $N_n\sim\dpois(n\nu_0)$, and $W_i\sim W$ are independent
random variables. Then
%
%e3.23 #&#
\begin{equation}
\label{ePoisson-bound} \mathbb{P}(\Omega_{n,x}) \le\int_{-\infty}^{x/2}
\mathbb{P}\{\zeta_{N_n} \in x - z + I\} \mathbb{P}\{V_n \in
\mathrm {d}z\} + \mathbb{P}\{V_n\ge x/2\}.
\end{equation}
Since $\mu:= \ln\mean[e^W]<\infty$, $\mathbb{P}\{V_n\ge x/2\} \le
\mean[e^{V_n - x/2}] = e^{n\mu- x/2}$. Therefore,
%
%e3.24 #&#
\begin{equation}
\label{ePoisson-V} \max_{n\le\rvf(\delta x)} \mathbb{P}\{V_n\ge x/2
\} \le\max_{n\le\rvf(\delta x)} e^{n\mu- x/2} = O \bigl(e^{-x/4}
\bigr).
\end{equation}
Next, for $z\le x/2$,
%
%e3.25 #&#
%e3.26 #&#
\begin{eqnarray}
\label{ePoisson-sum}&& \mathbb{P}\{\zeta_{N_n} \in x-z+I\}
\nonumber
\\[-8pt]
\\[-8pt]
\nonumber
&&\qquad\le\sum_{k>n\nu_0/2} \mathbb{P}\{\zeta_k \in
x-z+I\} \mathbb{P}\{N_n = k\} + \mathbb{P}\{N_n \le n
\nu_0/2\}.
\end{eqnarray}
Since $\mean[e^{-N_n}] = e^{n\nu_0(1/e-1)}$, by Markov's
inequality,
%
%e3.27 #&#
\begin{equation}
\label{ePoisson-N} \max_{n\ge L(x)} \mathbb{P}\{N_n \le n
\nu_0/2\} \le\max_{n\ge L(x)} e^{-n\nu_0(1/2-1/e)} \le
e^{-L(x)\nu_0/10}.
\end{equation}
On the other hand, note that for $t>1$, $\mathbb{P}\{Y\in t+I\} =
\nu(t+I)/\nu_0$ and $\overline{G}(t) = \overline{\nu}(t)/\nu_0$.
Then, as
$x-z\ge x/2$ and $Y_i>1$, for each $k>n\nu_0/2\ge L(x)\nu_0/2$, by
Lemma~\ref{lE-bound},
\[
\mathbb{P}\{\zeta_k\in x-z+I\} \le\frac{C_\nu k\overline{\nu}
(x-z)}{x-z}
\biggl[T_0 + \frac{\tilde\iof(2 (x-z), T_0)}{a_k} \biggr],
\]
where $C_\nu$ is a constant only depending on $\nu$. Then by
\eqref{ePoisson-bound}--\eqref{ePoisson-N}, letting $x_n = x-V_n$,
\[
\mathbb{P}(\Omega_{n,x}) \le C_\nu\mean \biggl[\mathbf{1}
\{V_n \le x/2\} \frac{N_n\overline{\nu}(x_n)}{x_n} \biggl[T_0 +
\frac
{\tilde\iof
(2 x_n, T_0)}{a_{N_n}} \biggr] \biggr] + \rx_n(x),
\]
where $\max_{L(x)\le n\le\rvf(\delta x)} \rx_n(x) = o(x^{-M})$ for
any $M>0$. Note that $N_n$ and $V_n$ are independent. Since
$\overline{\nu}(x)/x$ is regularly varying and decreasing, then
\begin{eqnarray*}
\mathbb{P}(\Omega_{n,x})&\le& C_\nu T_0 n\frac{\overline{\nu}(x)}{x}\\
&&{} +
C_\nu' \mean \biggl[\frac{N_n}{a_{N_n}} \biggr] \mean
\biggl[ \mathbf{1}\{x_n\ge x/2\}\frac{\overline{\nu}(x_n)
\tilde\iof(2 x_n, T_0)}{x_n} \biggr] +
\rx_n(x).
\end{eqnarray*}
Since
\[
\mean \biggl[\frac{N_n}{a_{N_n}} \mathbf{1}\{N_n<n
\nu_0/2 \mbox{ or } N_n>2 n\nu_0\} \biggr] =
o \bigl(e^{-c n} \bigr),\qquad n \to\infty,
\]
where $c>0$ is a constant, then by dominated convergence,
\[
\frac{a_n}{n} \mean \biggl[\frac{N_n}{a_{N_n}} \biggr] \sim\mean \biggl[
\frac{N_n/n}{a_{N_n}/a_n}\mathbf{1}\{n\nu_0/2 \le N_n \le2n
\nu_0\} \biggr] \sim\nu_0^{1-1/\alpha}.
\]
Consequently,
\[
\mathbb{P}(\Omega_{n,x})
\le C_\nu T_0 n\frac{\overline{\nu}(x)}{x} + C_\nu
\frac{n}{a_n} \mean \biggl[\mathbf{1}\{x_n\ge x/2\}
\frac{\overline
{\nu}(x_n)
\tilde\iof(2 x_n, T_0)}{x_n} \biggr] + \rx_n(x).
\]
Starting at this point, the treatment is very similar to that
following \eqref{eE-n}. First, by
\begin{eqnarray*}
&& \frac{x}{\rvf(x)} \sum_{L(x)\le n\le\rvf(\delta x)} \mathbb{P}(
\Omega_{n,x})
\\
&&\qquad\le\frac{C_\nu T_0}{\rvf(x)^2} \sum_{n\le\rvf(\delta x)} n +
C_\nu\tilde R(x) + \frac{x}{\rvf(x)} \sum
_{L(x)\le n\le\rvf(\delta x)} \rx_n(x)
\\
&&\qquad= O \bigl(\delta^2 \bigr) + C_\nu\tilde R(x) + o(1),
\end{eqnarray*}
where writing $\theta=1/\alpha-1$,
\[
\tilde R(x) = \frac{x}{\rvf(x)} \sum_{L(x)\le n\le\rvf(\delta x)}
\frac{n^{-\theta}}{\svf(n)} \mean \biggl[\mathbf{1}\{x_n\ge x/2\}
\frac{\overline{\nu}(x_n)
\tilde\iof(2 x_n, T_0)}{x_n} \biggr].
\]
If $\alpha\in(0,1/2)$, then by Assumption~\ref{abasic},
\eqref{ediff2} holds for $\tilde\iof(x,T_0)$, and hence
\begin{eqnarray*}
\tilde R(x) &=& \frac{x}{\rvf(x)} o \biggl( \sum_{L(x)\le n\le\rvf
(\delta x)}
\frac{n^{-\theta}}{\svf(n)} \mean \biggl[\mathbf{1}\{x_n\ge x/2\}
\frac{\overline{\nu}(x_n)
\rvf(2 x_n)^2}{x_n u_\theta(2 x_n)} \biggr] \biggr)
\\
&=& \frac{x}{\rvf(x)} o \biggl( \sum_{L(x)\le n\le\rvf(\delta x)}
\frac{n^{-\theta}}{\svf(n)} \frac{\rvf(x)}{x u_\theta(x)} \biggr) = o(1),\qquad x\to\infty.
\end{eqnarray*}
If $\alpha=1/2$, then by Assumption~\ref{abasic}, \eqref{ediff}
holds for $\tilde\iof(x,T_0)$. As a result, if $\tilde u(x)/\tilde
u(L(x))\to1$, then
\begin{eqnarray*}
\tilde R(x) &=& O \biggl( \sum_{L(x)\le n\le\rvf(\delta x)}
\frac
{n^{-\theta}}{\svf(n)} \mean \biggl[\mathbf{1}\{x_n\ge x/2\}
\frac{x_n}{\rvf(x_n)} \cdot \frac{\overline{\nu}(x_n)\rvf(2 x_n)^2}{x_n \tilde u(2 x_n)} \biggr] \biggr)
\\
&=& O \biggl( \sum_{L(x)\le n\le\rvf(\delta x)} \frac{n^{-\theta
}}{\svf(n)}
\frac{1}{\tilde
u(x)} \biggr) = o(1),\qquad x\to\infty.
\end{eqnarray*}
The case $\tilde u(x)/\tilde u(L(x))\not\to1$ can be shown likewise.
This then completes the proof of Theorem~\ref{trenewal-ind}.

\begin{pf*}{Proof of Lemma~\ref{lind}}
By Proposition~\ref{pprior}, it suffices to show that as $x\to\infty $,
$\mathbb{P}\{X\in x + I\} = O(\overline{\cdf}(x)/x^{1-\kappa})$.
For $x>0$,
\begin{eqnarray*}
&& \mathbb{P}\{X\in x + I\}
\\
&&\qquad= \mathbb{P}\{\zeta_N + W \in x + I\}
\\
&&\qquad \le\mathbb{P}\{\zeta_N + W\in x + I, N<\ln x, W<x/2\} + \mathbb
{P}\{W\ge x/2\}
\\
&&\qquad\quad{} + \mathbb{P}\{N\ge\ln x\}
\\
&&\qquad \le\sum_{n<\ln x} \frac{e^{-\nu_0} \nu_0^n}{n!} \sup
_{z>x/2} \mathbb{P}\{\zeta_n \in z + I\} + \mathbb{P}
\{W\ge x/2\} + \mathbb{P}\{N\ge\ln x\}.
\end{eqnarray*}
Given $\gamma\in(0,1)$, by $x^\gamma\ln x = o(x)$, for $x\gg1$,
$n<\ln x$ and $z>x/2$, if $\zeta_n\in z+I$, then there is at
least one $1\le i\le n$ with $Y_i > z^\gamma$, and if there is
exactly one such $i$, then $Y_i>z/2$. Thus
\begin{eqnarray*}
&& \mathbb{P}\{\zeta_n \in z + I\}
\\
&&\qquad\le n^2 \mathbb{P}\bigl\{\zeta_n\in z+I,
Y_{n-1}>z^\gamma, Y_n > z^\gamma\bigr\} +
n\mathbb{P}\{\zeta_n\in z+I, Y_n>z/2\}.
\end{eqnarray*}
First, following the argument to bound $E\sp{(3)}_{n,x}$ in the proof of
Lemma~\ref{lE-bound},
\begin{eqnarray*}
&& \mathbb{P}\bigl\{\zeta_n\in z+I, Y_{n-1}>z^\gamma,
Y_n > z^\gamma\bigr\}
\\
&&\qquad= \sum_{k=0}^\infty \mathbb{P}\bigl\{
\zeta_n\in z+I, Y_{n-1}>z^\gamma, Y_n
\in z^\gamma+ k h +I\bigr\}
\\
&&\qquad\le\sup_{t>z^\gamma} \mathbb{P}\{Y\in t +I\} \sum
_{k=0}^\infty \mathbb{P}\bigl\{ \zeta_{n-1}
\in z- z^\gamma- k h +J, Y_{n-1}>z^\gamma\bigr\}
\\
&&\qquad\le2\mathbb{P}\bigl\{Y>z^\gamma\bigr\} \sup_{t>z^\gamma}
\mathbb{P}\{Y \in t +I\}.
\end{eqnarray*}
By $\mathbb{P}\{Y>z^\gamma\} = \overline{\nu}(z^\gamma)/\nu_0$ and
$\mathbb{P}\{Y \in t +
I\} = \nu(t+I)/\nu_0 = O(\overline{\nu}(t)/t^{1-\kappa})$, the RHS
of the display is $O(\overline{\nu}(z^\gamma)^2 / z^{\gamma
(1-\kappa)})
= O(\overline{\nu}(x^\gamma)^2 / x^{\gamma(1-\kappa)})$. It
follows that
if $\gamma>(\alpha+ 1-\kappa)/(2\alpha+ 1-\kappa)$, then the RHS
is $o(\overline{\cdf}(x)/ x^{1-\kappa})$. With a similar argument, we
also get
\[
\mathbb{P}\{\zeta_n\in z+I, Y_n>z/2\} \le2 \sup
_{t>z/2}\mathbb{P}\{Y \in z/2+I\} = O \bigl(\overline{\cdf
}(x)/x^{1-\kappa} \bigr).
\]
As a result,
\begin{eqnarray*}
\sum_{n<\ln x} \frac{e^{-\nu_0} \nu_0^n}{n!} \sup
_{z>x/2} \mathbb{P}\{\zeta_n \in z+I\} &=& O \bigl(
\mean N^2\cdot\overline{\cdf}(x)/x^{1-\kappa} \bigr)
\\
&=& O \bigl(\overline{\cdf}(x)/x^{1-\kappa} \bigr).
\end{eqnarray*}
On the other hand, $\mathbb{P}\{W\ge x/2\} \le\mean[e^{2(W-x/2)}] =
O(e^{-x})$ and, for any $M>0$, $\mathbb{P}\{N\ge\ln x\} \le\mean
[e^{M(N-\ln
x)}] = O(x^{-M})$. By letting $M>\alpha+ {1-\kappa}$, the above
bounds together yield $\mathbb{P}\{X\in x+I\} =
O(\overline{\cdf}(x)/x^{1-\kappa})$, as desired.
\end{pf*}

%
%sA #&#
\begin{appendix}\label{app}
\section*{Appendix}
\setcounter{equation}{0}
\renewcommand{\thetheorem}{A}
In Section~\ref{sresults}, we remarked that if the SRT holds, then
\eqref{elow-cut} holds for any $L(x)=o(\rvf(x))$. This follows from
the following

%
%prA.1 #&#
\begin{prop} \label{pweak-renewal}
For $\cdf$ satisfying both \eqref{etail-svf} and
\eqref{etail-ratio},
%
%eA.1 #&#
\begin{equation}
\label{erenewal-inf} \liminf_{x\to\infty } x\overline{\cdf}(x) \renf(x+I) =
h \Lambda_\cdf,
\end{equation}
where $\Lambda_\cdf$ is defined in \eqref{erenewal}, and $h>0$ is
arbitrary if $\cdf$ is nonarithmetic and is the span of $\cdf$
otherwise.
\end{prop}

Indeed, if the SRT holds, then $\liminf$ in \eqref{erenewal-inf} can
be replaced with $\lim$. On the other hand, by Lemma~\ref{lmiddle-n},
%
%eA.2 #&#
\begin{equation}
\label{etruncate} \lim_{\delta\to0} \lim_{x\to\infty } x
\overline{\cdf}(x) \sum_{n\in J_\delta(x)} \mathbb{P}
\{S_n\in x + I\} = h\Lambda_\cdf,
\end{equation}
where $J_\delta(x) = (\rvf(\delta x), \rvf(x/\delta))$. It then
follows that
\[
\lim_{\delta\to0} \lim_{x\to\infty } x\overline{\cdf}(x)
\sum_{n\le\rvf(\delta x)} \mathbb{P}\{S_n\in x + I\} = 0
\]
and hence \eqref{elow-cut} holds for any $L(x) = o(\rvf(x))$.

\begin{pf*}{Proof of Proposition~\ref{pweak-renewal}}
It is well known that if $\cdf$ is nonlattice with support in
$[0,\infty)$ and infinite mean, then Proposition~\ref
{pweak-renewal} holds (\cite{bingham89}, Theorem~8.6.6). For the
general case, we follow the proof in \cite{bingham89}. Denoting
$V(x) = \renf((0, x^+])$, the starting point is the identity
%
%eA.3 #&#
\begin{equation}
\label{eweak-renewal} \lim_{x\to\infty } \overline{\cdf}(x) V(x) =
\Lambda_\cdf/\alpha.
\end{equation}
This is established on page~361 in \cite{bingham89}. However, the
proof there relies on the Laplace transforms of $\cdf$ and $\renf$,
so it cannot apply to the general case as the transforms may be
$\infty$. Instead, we shall prove \eqref{eweak-renewal} using a
more probabilistic argument, which is basically a coarse version of
the one for the SRT. For now assume \eqref{eweak-renewal} to be
true. Then \eqref{etruncate} implies
\[
\liminf_{x\to\infty } x\overline{\cdf}(x) \renf(x+I) \ge h
\Lambda_\cdf.
\]
Assume that strict inequality holds. Then by \eqref
{eweak-renewal}, there is $h'>h$, such that for all $x\gg1$,
$\renf(x+I) \ge h'\alpha V(x)/x$. Also $V\in\CR_\alpha$. Then
as $t\to\infty $,
\[
\int_0^t \renf(x+I) \,\mathrm{d}x \ge
\bigl(1+o(1) \bigr) h'\alpha\int_0^t
V(x) x^{-1} \,\mathrm{d}x \sim h' V(t).
\]
However, since $\renf(x+I) = V(x+h) - V(x)$ for $x\ge0$, LHS
$= \int_t^{t+h} V - \int_0^h V \sim h V(t)$, which contradicts the
above display. Thus \eqref{erenewal-inf} follows.

It remains to show \eqref{eweak-renewal}. Given $\delta\in(0,1)$,
%
%eA.4 #&#
\begin{equation}
\label{etruncate-low} \sum_{n\le\rvf(\delta x)} \mathbb{P}\bigl
\{S_n \in(0,x]\bigr\} \le\rvf(\delta x)\sim\delta^\alpha
\rvf(x),\qquad x\to \infty.
\end{equation}
On the other hand, by the LLTs and the boundedness of $g$, there is
$C>0$, such that for all $x\gg1$, $n\ge
\rvf(x/\delta)$, and $t\in\Reals$, $\mathbb{P}\{S_n \in t+I\} \le C/a_n$.
Then, by dividing $(0,x]$ into $\lceil x/h\rceil$ intervals of equal
length, it is seen that $\mathbb{P}\{S_n\in(0,x]\} \le C x/a_n$.
Consequently,
%
%eA.5 #&#
\begin{eqnarray}
\label{etruncate-up} \sum_{n\ge\rvf(x/\delta)} \mathbb{P}\bigl
\{S_n \in(0,x]\bigr\} &\le &C x\sum_{n\ge\rvf(x/\delta)}
\frac{1}{a_n}
\nonumber
\\[-8pt]
\\[-8pt]
\nonumber
& \sim&\frac{C'x\rvf(x/\delta)}{x/\delta} \sim C'''
\delta^{1-\alpha}\rvf(x).
\end{eqnarray}

By the central limit theorem, as $n\to\infty $, $G_n(s):= \mathbb
{P}\{0<
S_n/a_n\le s\}\to G(s):= \int_0^{s^+} g$ for each $s$. Since $G_n$
and $G$ are nondecreasing functions with range contained in
$[0,1]$, and $G$ is continuous, the pointwise convergence gives
$\sup|G_n - G|\to0$. Then by $\mathbb{P}\{S_n\in(0,x]\} =
G_n(x/a_n)$, as
$x\to\infty $,
\begin{eqnarray*}
&&\sum_{n\in J_\delta(x)} \mathbb{P}\bigl\{S_n
\in(0,x]\bigr\}\\
&&\qquad= \sum_{n \in J_\delta(x)} G(x/a_n) +
\bigl[\rvf(x/\delta) - \rvf(\delta x) \bigr]o(1)
\\
&&\qquad= \int_{\rvf(\delta x)}^{\rvf(x/\delta)} G \bigl(x/
\rvf^{-1}(t) \bigr) \,\mathrm{d}t + \bigl(\delta^{-\alpha} -
\delta^\alpha \bigr) o(1) \rvf(x).
\end{eqnarray*}
By change of variable $u = x/\rvf^{-1}(t)$ and $\rvf'(u) \sim
\alpha\rvf(u)/u$ as $u\to\infty $,
\begin{eqnarray*}
\int_{\rvf(\delta x)}^{\rvf(x/\delta)} G \bigl(x/\rvf^{-1}(t)
\bigr) \,\mathrm{d}t &=& \int_\delta^{1/\delta} G(u)
\frac{x}{u^2} \rvf'(x/u) \,\mathrm{d}u
\\
&\sim&\int_\delta^{1/\delta} G(u)\frac{\alpha\rvf(x/u)}{u}
\,\mathrm{d}u
\\
&\sim&\rvf(x)\int_\delta^{1/\delta} G(u)\alpha
u^{-1-\alpha} \,\mathrm{d}u,\qquad x\to\infty.
\end{eqnarray*}
As a result,
%
%eA.6 #&#
\begin{equation}
\label{etruncate-mid} \lim_{\delta\to0} \lim_{x\to\infty }
\frac{1}{
\rvf(x)} \sum_{n\in J_\delta(x)} \mathbb{P}\bigl
\{S_n\in(0,x]\bigr\} = \int_0^\infty
G(u)\alpha u^{-1-\alpha} \,\mathrm{d}u.
\end{equation}
The RHS is $\Lambda_\cdf/\alpha$. Combining
\eqref{etruncate-low}--\eqref{etruncate-mid}, then
\eqref{eweak-renewal} follows.
\end{pf*}
\end{appendix}

\section*{Acknowledgments} The author would like to thank
the referees for careful reviews and useful suggestions.

%
%
% imsref loaded by akundreckaite, 2014-05-26 14:47:16
%

%


%suskaldyti doi

\printaddresses
\end{document}